\documentstyle[amssymb,amsfonts]{amsart}

\def\R{{\hbox{\bf R}}}
\def\C{{\hbox{\bf C}}}

\def\B{{\hbox{\bf B}}}
\def\T{{\hbox{\bf T}}}
\def\E{{\hbox{\bf E}}}
\def\ex{{\mathcal{E}}}
\def\re{{\mathcal{R}}}
\def\ka{{\mathcal{K}}}
\def\tubes{{\mathbb{T}}}\def\vi{{\mathbb{V}}}
\def\v{{\hbox{\bf v}}}
\def\w{{\hbox{\bf w}}}
\def\A{{\hbox{\bf A}}}

\font \roman = cmr10 at 10 true pt

\def\supp{{\operatorname{supp}}}

\def\im{{\hbox{\roman im}}}

\def\adj{{\hbox{\roman adj}}}

\def\Z{{\hbox{\bf Z}}}
\def\N{{\hbox{\bf N}}}

\def\th{{\operatorname{th}}}
\def\eps{\varepsilon}
\newenvironment{proof}{\noindent {\bf Proof} }{\endprf\par}
\def \endprf{\hfill  {\vrule height6pt width6pt depth0pt}\medskip}
\def\emph#1{{\it #1}}
\def\textbf#1{{\bf #1}}

\theoremstyle{plain}
  \newtheorem{theorem}[subsection]{Theorem}
  \newtheorem{conjecture}[subsection]{Conjecture}

  \newtheorem{proposition}[subsection]{Proposition}
  
  \newtheorem{lemma}[subsection]{Lemma}
  \newtheorem{corollary}[subsection]{Corollary}
  \newtheorem{question}[subsection]{Question}
\theoremstyle{remark}
  \newtheorem{remark}[subsection]{Remark}

\theoremstyle{definition}

\include{psfig}
\thanks{The first author was supported by an EPSRC Postdoctoral Fellowship,
the second by EC project ``HARP'', and the third by a grant from the Packard Foundation}
\begin{document}
\title{On the Multilinear Restriction and Kakeya conjectures}
\date{August 2005}
\begin{abstract}
We prove $d$-linear analogues of the
classical restriction and Kakeya conjectures in $\R^d$. Our
approach involves obtaining monotonicity formulae pertaining to a
certain evolution of families of gaussians, closely related to
heat flow. We conclude by giving some applications to the
corresponding variable-coefficient problems and the so-called
``joints" problem, as well as presenting some $n$-linear analogues for $n<d$.
\end{abstract}

\author{Jonathan Bennett}
\address{School of Mathematics, University of Birmingham, Birmingham, B15 2TT, UK}
\email{J.Bennett@@bham.ac.uk}

\author{Anthony Carbery}
\address{School of Mathematics, University of Edinburgh, Edinburgh, EH9 3JZ, UK}
\email{A.Carbery@@ed.ac.uk}

\author{Terence Tao}
\address{Department of Mathematics, UCLA, Los Angeles CA 90095-1555}
\email{tao@@math.ucla.edu}

\maketitle
\section{Introduction}
For $d\geq 2$, let $U$ be a compact neighbourhood of the origin in
$\R^{d-1}$ and $\Sigma :U\rightarrow \R^{d}$ be a smooth
parametrisation of a $(d-1)$-dimensional submanifold $S$ of
$\R^d$ (for instance, $S$ could be a small portion of the unit sphere $S^{d-1}$).
To $\Sigma$ we associate the \textit{extension operator}
$\ex$, given by
$$
\ex g(\xi):=\int_{U}g(x)e^{i\xi\cdot\Sigma(x)}dx,
$$
where $g\in L^{1}(U)$ and $\xi\in\R^d$. This operator is sometimes
referred to as the \textit{adjoint restriction operator} since its
adjoint $\ex^{*}$ is given by $\ex^{*}f=\widehat{f}\circ\Sigma$,
where $\widehat{\;}\:$ denotes the $d$-dimensional Fourier
transform. It was observed by E. M. Stein in the late 1960's that
if the submanifold parametrised by $\Sigma$ has everywhere
non-vanishing gaussian curvature, then non-trivial
$L^{p}(U) \to L^{q}(\R^{d})$ estimates for $\ex$ may be obtained.
The classical \textit{restriction conjecture} concerns the full
range of exponents $p$ and $q$ for which such bounds hold.
\begin{conjecture}[Linear Restriction]\label{LRC}
If $S$ has everywhere non--vanishing gaussian curvature,
$q>\tfrac{2d}{d-1}$ and $p'
\leq\tfrac{d-1}{d+1}q$, then there exists a constant
$0 < C < \infty$ depending only on $d$ and $\Sigma$ such that
\begin{equation*}
\|\ex g\|_{L^{q}({\bf R}^d)}\leq C\|g\|_{L^{p}(U)}
\end{equation*}
for all $g\in L^p(U)$.
\end{conjecture}
See for
example \cite{TaoRestn} for a discussion of the progress made on
this problem, the rich variety of techniques that have
developed in its wake, and the connection to other problems in harmonic analysis,
partial differential equations, and geometric analysis.  In particular the restriction
problem is intimately connected to the \emph{Kakeya problem}, which we shall discuss later in this
introduction.

In recent years certain \textit{bilinear} analogues of the
restriction problem have come to light in natural ways from a
number of sources (see for example \cite{Bolattice},
\cite{KM}, \cite{MVV},
\cite{TVV}, \cite{wolffcone}, \cite{Taocone}, \cite{TaoKM}, \cite{vargas},
\cite{Lee},
\cite{Erd} concerning the
well-posedness theory of non-linear dispersive
equations, and applications to a variety of central problems in harmonic and
geometric analysis). More
specifically, given two such smooth mappings $\Sigma_{1}$ and
$\Sigma_{2}$, with associated extension operators $\ex_{1}$ and
$\ex_{2}$, one may ask for which values of the exponents $p$ and
$q$, the bilinear operator
$(g_{1},g_{2})\mapsto\ex_{1}g_{1}\:\ex_{2}g_{2}$ may be bounded
from $L^{p}\times L^{p}$ to $L^{q/2}$. The essential point here is
that if the submanifolds parametrised by $\Sigma_{1}$ and
$\Sigma_{2}$ are assumed to be \textit{transversal} (up to
translations), then one can expect the range of such exponents to
broaden; again see \cite{TVV}. However, one of the more puzzling
features of such bilinear problems is that, in three dimensions
and above, they seem to somewhat confuse the role played by the
curvature of the associated submanifolds. For example, it is known
that the bilinear restriction theories for the cone and paraboloid
are almost identical, whereas the linear theories for these surfaces
are not (see \cite{TaoRestn} for further discussion of this).
Moreover, simple heuristics suggest that the optimal ``$k$-linear"
restriction theory requires at least $d-k$ nonvanishing principal
curvatures, but that further curvature assumptions have no further effect.
In $d$ dimensions it thus seems particularly natural
to consider a $d$-linear set-up, as one then does not expect to require any curvature conditions.

For each $1\leq j\leq d$ let $\Sigma_{j}:U_{j}\rightarrow \R^{d}$
be such a smooth mapping and let $\ex_j$ be the associated
extension operator. Our analogue of the bilinear transversality
condition will essentially amount to requiring that the normals to
the submanifolds parametrised by the $\Sigma_{j}$'s span at all
points of the parameter space. In order to express this in an
appropriately uniform manner let $A,\nu>0$ be given, and for each
$1\leq j\leq d$ let $Y_j$ be the $(d-1)$-form
$$
Y_{j}(x):=\bigwedge_{k=1}^{d-1}\frac{\partial}{\partial
x_{k}}\Sigma_{j}(x)
$$
for all $x \in U_j$; by duality we can view $Y_j$ as a vector field on $U_j$.
We will not impose any curvature conditions (in particular, we permit the vector fields
$Y_j$ to be constant), but we will impose the ``transversality" (or ``spanning") condition
\begin{equation}\label{spanning}
\det\left(Y_{1}(x^{(1)}),\hdots,Y_{d}(x^{(d)})\right)\geq\nu,
\end{equation}
for all $x^{(1)}\in U_{1},\hdots,x^{(d)}\in U_d$, along with the
smoothness condition
\begin{equation}\label{smoothcond}
\|\Sigma_{j}\|_{C^{2}(U_j)}\leq A\;\mbox{ for all }\;1\leq j\leq d.
\end{equation}

\begin{remark}\label{alt}
If $U_j$ is sufficiently small then
$\ex_{j}g_{j}=\widehat{G_{j}d\sigma_{j}}$, where
$G_{j}:\Sigma_{j}(U_{j})\rightarrow\C$ is the ``normalised lift"
of $g_{j}$, given by
$G_{j}(\Sigma_{j}(x))=|Y_{j}(x)|^{-1}g_{j}(x)$, and $d\sigma_{j}$
is the induced Lebesgue measure on $\Sigma_{j}(U_{j})$.
\end{remark}

By testing on the standard examples that generate the
original linear restriction conjecture (characteristic functions of
small balls in $\R^{d-1}$ -- see \cite{St2}) we are led to the
following conjecture\footnote{Strictly speaking, this is a ``multilinear extension'' or ``multilinear
adjoint restriction'' conjecture rather than a multilinear restriction conjecture, but the use of the term ``restriction'' is
well established in the literature.}.
\begin{conjecture}[Multilinear Restriction]\label{MLRC}
Suppose that \eqref{spanning} and \eqref{smoothcond} hold,
$q\geq\tfrac{2d}{d-1}$ and $p'\leq \tfrac{d-1}{d}q$. Then there
exists a constant $C$, depending only on $A$, $\nu$, $d$, and $U_1,\ldots,U_d$, for
which
\begin{equation}\label{mr}
\Bigl\|\prod_{j=1}^{d}\ex_{j}g_{j}\Bigr\|_{L^{q/d}({\bf R}^d)}\leq
C\prod_{j=1}^{d} \|g_{j}\|_{L^{p}(U_{j})}
\end{equation}
for all $g_{1}\in L^p(U_1),\hdots,g_{d}\in L^{p}(U_d)$.
\end{conjecture}
\begin{remark}\label{scaling}
Using a partition of unity and an appropriate affine
transformation we may assume that $\nu\sim 1$ and that for each
$1\leq j\leq d$, $\Sigma_{j}(U_{j})$ is contained in a
sufficiently small neighbourhood of the $j^\th$ standard basis
vector $e_{j}\in\R^d$.
\end{remark}
\begin{remark}\label{Loomis}
By multilinear interpolation (see for example \cite{bergh:interp})
and H\"older's inequality, Conjecture \ref{MLRC}
may be reduced to the endpoint case $p=2$, $q=\tfrac{2d}{d-1}$;
i.e. the $L^{2}$ estimate
\begin{equation}\label{endpoint}
\Bigl\|\prod_{j=1}^{d}\ex_{j}g_{j}\Bigr\|_{L^{2/(d-1)}({\bf
R}^d)}\leq C\prod_{j=1}^{d} \|g_{j}\|_{L^{2}(U_j)}.
\end{equation}
We emphasise that at this $d$-linear level, the optimal estimate is on
$L^2$, rather than $L^{\frac{2d}{d-1}}$.
It should also be pointed out that the conjectured range of
exponents $p$ and $q$ is independent of any additional curvature
assumptions that one might make on the submanifolds parametrised
by the $\Sigma_{j}$'s. This is very much in contrast with similar
claims at lower levels of multilinearity. It is instructive to
observe that if the mappings $\Sigma_{j}$ are \textit{linear},
then by an application of Plancherel's theorem, the conjectured
inequality \eqref{endpoint} (for an appropriate constant $C$)
is equivalent to the
classical Loomis--Whitney inequality \cite{LW}. This elementary
inequality states that if $\pi_{j}:\R^{d}\rightarrow\R^{d-1}$ is
given by $\pi_{j}(x):=(x_{1},\hdots,x_{j-1},x_{j+1},\hdots,x_{d})$,
then
\begin{eqnarray}\label{clw}
\begin{aligned}
\int_{{\bf R}^{d}}f_{1}(\pi_{1}(x))\cdot\cdot\cdot
f_{d}(\pi_{d}(x)) \;dx \leq \|f_{1}\|_{d-1}\cdot\cdot\cdot
\|f_{d}\|_{d-1}
\end{aligned}
\end{eqnarray}
for all $f_{j}\in L^{d-1}(\R^{d-1})$. One may therefore view the
multilinear restriction conjecture as a certain (rather oscillatory)
generalisation of the Loomis--Whitney inequality.
The nature of this generalisation is
clarified in Section \ref{RK}.
\end{remark}
\begin{remark}
Conjecture \ref{MLRC} in two dimensions is elementary and
classical, and is implicit in arguments of C. Fefferman and Sj\"olin.
In three dimensions this
(trilinear) problem was considered in \cite{BCW} (and previously
in \cite{B}), where some partial results on the sharp line $p'=
\tfrac{d-1}{d}q$ were obtained.
\end{remark}

It is a well-known fact that the linear
restriction conjecture implies the so-called (linear)
\emph{Kakeya conjecture}.  This conjecture takes several forms.
One particularly simple one is the assertion that any (Borel) set in $\R^n$
which contains a unit line segment in every direction must have full Hausdorff (and
thus Minkowski) dimension.  Here we shall consider a more quantitative version
of the conjecture, which is stronger than the one just described.
For $0<\delta\ll 1$ we define a
$\delta$-\textit{tube} to be any rectangular box $T$ in $\R^{d}$ with $d-1$
sides of length $\delta$ and one side of length $1$; observe that such tubes have
volume $|T| \sim \delta^{d-1}$. Let
$\mathbb{T}$ be an arbitrary collection of such $\delta$-tubes
whose orientations form a $\delta$-separated set of points on
${\bf S}^{d-1}$.
We use $\# \mathbb{T}$ to denote the cardinality of $\mathbb{T}$, and $\chi_T$ to denote the indicator
function of $T$ (thus $\chi_T(x) = 1$ when $x \in T$ and $\chi_T(x)=0$ otherwise).

\begin{conjecture}[Linear Kakeya]\label{LKC}  Let $\mathbb{T}$ and $\delta$ be as above.
For each $\tfrac{d}{d-1}<q\leq\infty$ there is a constant $C$,
independent of $\delta$ and the collection $\mathbb{T}$, such that
\begin{equation}\label{kakeya}
\Bigl\|\sum_{T\in\mathbb{T}}\chi_{T}\Bigr\|_{L^{q}({\bf
R}^{d})} \leq C\delta^{(d-1)/q}\left(
\# \mathbb{T}\right)^{1/q}.
\end{equation}
\end{conjecture}
The proof that Conjecture \ref{LRC} implies Conjecture \ref{LKC}
follows a standard
Rademacher-function argument going back implicitly to \cite{Feff}
and \cite{BCSS}.
The endpoint $q=\tfrac{d}{d-1}$ of \eqref{kakeya} can be seen to be
false (unless one places an additional
logarithmic factor in $\delta$ on the right-hand side), either
by considering a collection of tubes passing through
the origin, or by Besicovitch set examples.
See \cite{Wolffbook} for a detailed account of these facts.

By a straightforward adaptation of the techniques in the linear
situation, the multilinear restriction conjecture can be seen to
imply a corresponding multilinear Kakeya-type conjecture. Suppose
$\tubes_{1},\hdots,\tubes_{d}$ are families of $\delta$-tubes in
$\R^d$. We allow the tubes within a single family $\tubes_j$ to be parallel.
However, we assume that for each $1\leq j\leq d$, the tubes
in $\tubes_{j}$ have long sides pointing in directions belonging
to some sufficiently small \textit{fixed} neighbourhood of the
$j^{\th}$ standard basis vector $e_{j}$ in ${\bf S}^{d-1}$. It will
be convenient to refer to such a family of tubes as being
\textit{transversal}. (The vectors $e_{1},\hdots,e_{d}$ may be
replaced by any fixed linearly independent set of vectors in
$\R^d$ here, as affine invariance considerations reveal.)
\begin{conjecture}[Multilinear Kakeya]\label{MLKC}
Let $\T_1,\ldots,\T_d$ and $\delta$ be as above.
If $\tfrac{d}{d-1}\leq q\leq\infty$ then there exists a constant
$C$, independent of $\delta$ and the families of tubes
$\tubes_{1},\hdots,\tubes_{d}$, such that
\begin{equation}\label{mk}
\Bigl\|\prod_{j=1}^{d}\Bigl(\sum_{T_{j}\in\tubes_{j}}\chi_{T_{j}}\Bigr)\Bigr\|_{L^{q/d}({\bf
R}^d)}\leq C\prod_{j=1}^{d}(\delta^{d/q}\: \mbox{\#}\tubes_{j}).
\end{equation}
\end{conjecture}
\begin{remark}
Since the case $q=\infty$ is trivially true, the above conjecture
is equivalent via H\"older's inequality
to the endpoint case $q=\tfrac{d}{d-1}$.  In contrast to the linear setting,
there is no obvious counterexample prohibiting this claim holding at the endpoint
$q = \tfrac{d}{d-1}$, and indeed in the $d=2$ case it is easy to verify this endpoint estimate.
\end{remark}
\begin{remark}\label{superficial}
By contrast with similar statements at lower levels of multilinearity,
each family $\tubes_j$ is permitted to contain parallel tubes, and
even arbitrary repetitions of tubes. By scaling and a limiting argument
we thus see that
the conjectured inequality reduces
to the superficially
stronger
\begin{equation}\label{coefficients}
\Bigl\|\prod_{j=1}^{d}\Bigl(\sum_{T_{j}\in\tubes_{j}}\chi_{T_{j}}*\mu_{T_{j}}\Bigr)\Bigr\|_{L^{q/d}({\bf
R}^d)} \leq C \prod_{j=1}^{d}\Bigl(\delta^{d/q}\:
\sum_{T_{j}\in\tubes_{j}}\|\mu_{T_{j}}\|\Bigr)
\end{equation}
for all finite measures $\mu_{T_{j}}$ ($T_j\in\tubes_j$,
$1\leq j\leq d$) on $\R^d$.
\end{remark}
\begin{remark}
The decision to formulate Conjecture \ref{MLKC} in term of
$\delta\times\cdots\times\delta\times 1$ tubes is largely for
historical reasons. However, just by scaling, it is easily seen
that \eqref{mk} is equivalent to the inequality
\begin{equation}\label{1tubes}
\Bigl\|\prod_{j=1}^{d}\Bigl(\sum_{\widetilde{T}_{j}\in\widetilde{\tubes}_{j}}
\chi_{\widetilde{T}_{j}}\Bigr)\Bigr\|_{L^{q/d}({\bf
R}^d)}\leq C\prod_{j=1}^{d}(\mbox{\#}\widetilde{\tubes}_{j}),
\end{equation}
where the collections $\widetilde{\tubes}_{j}$ consist of tubes of
width $1$ and arbitrary (possibly infinite) length. (Of course we
continue to impose the appropriate transversality condition on the
families $\widetilde{\tubes}_{1},\hdots,\widetilde{\tubes}_{d}$
here.)
\end{remark}

\begin{remark}
As may be expected given Remark \ref{Loomis}, the special case of
the conjectured inequality (or rather the equivalent form
\eqref{coefficients} with $q=\tfrac{d}{d-1}$, and an appropriate
constant $C$) where all of the
tubes in each family $\tubes_j$ are \textit{parallel}, is easily
seen to be equivalent to the
Loomis--Whitney inequality. We may therefore also view the
multilinear Kakeya conjecture as a generalisation of the
Loomis--Whitney inequality. The geometric nature of this
generalisation is of course much more transparent than that of
Conjecture \ref{MLRC}. In particular, one may find it enlightening
to reformulate \eqref{coefficients} (with $q=\tfrac{d}{d-1}$) as
an $\ell^1$ vector-valued version of \eqref{clw}.
\end{remark}

\begin{remark} As mentioned earlier, the linear Kakeya conjecture implies
something about the dimension of sets which contain a unit line segment in every
direction.  The multilinear Kakeya conjecture does not have a similarly simple geometric
implication, however there is a connection in a similar spirit between this conjecture and the joints problem;
see Section \ref{joints-sec}.
\end{remark}

Remarkably, at this $d$-linear level it turns out that the
restriction and Kakeya conjectures are essentially
\textit{equivalent}. This ``equivalence", which is the subject of
Section 2, follows from multilinearising a well-known
induction-on-scales argument of Bourgain \cite{Bo} (see also \cite{TVV}
for this argument in the bilinear setting). Once we have this
equivalence we may of course focus our attention on Conjecture
\ref{MLKC}, the analysis of which is the main innovation of this
paper. The general idea behind our approach to this conjecture is
sufficiently simple to warrant discussion here in the
introduction. First let us observe that if each $T_j\in\tubes_j$
is centred at the origin (for all $1\leq j\leq d$), then the
left and right hand sides of the conjectured
inequality \eqref{mk} are trivially comparable. This observation
leads to the suggestion that such
configurations of tubes might actually be \textit{extremal} for
the left hand side of \eqref{mk}.
\begin{question}
Is it reasonable to expect a quantity such as
$$
\Bigl\|\prod_{j=1}^{d}\Bigl(\sum_{T_{j}\in\tubes_{j}}\chi_{T_{j}}\Bigr)\Bigr\|_{L^{q/d}({\bf
R}^d)}$$
to be monotone increasing for $q\geq\tfrac{d}{d-1}$ as the
constituent tubes ``simultaneously slide" to the origin?
\end{question}
For reasons both analytic and algebraic, in pursuing this idea it
seems natural to replace the rough characteristic functions of
tubes by gaussians (of the form $e^{-\pi\langle
A(x-v),(x-v)\rangle}$ for appropriate positive definite $d\times
d$ matrices $A$ and vectors $v\in\R^d$) adapted to them. As we
shall see in Sections 3 and 4, with this gaussian reformulation
the answer to the above question is, to all intents and purposes,
yes for $q>\tfrac{d}{d-1}$.
In Section 3 we illustrate this by giving a new proof of the
Loomis--Whitney inequality, which we then are able to perturb in
Section 4.
As a corollary of our perturbed result in Section 4 we
obtain the multilinear Kakeya conjecture up to the endpoint, and a
``weak" form of the multilinear restriction conjecture.


More precisely, our main results are as follows.

\vspace{3mm}

\begin{theorem}[Near-optimal multilinear Kakeya]\label{kakthm}
If $\tfrac{d}{d-1}< q\leq\infty$ then there exists a constant
$C$, independent of $\delta$ and the transversal families of tubes
$\tubes_{1},\hdots,\tubes_{d}$, such that
\begin{equation*}
\Bigl\|\prod_{j=1}^{d}\Bigl(\sum_{T_{j}\in\tubes_{j}}\chi_{T_{j}}\Bigr)\Bigr\|_{L^{q/d}({\bf
R}^d)}\leq C\prod_{j=1}^{d}(\delta^{d/q}\: \mbox{\#}\tubes_{j}).
\end{equation*}
Furthermore, for each $\epsilon>0$ there is a similarly uniform constant $C$
for which
\begin{equation*}
\Bigl\|\prod_{j=1}^{d}\Bigl(\sum_{T_{j}\in\tubes_{j}}\chi_{T_{j}}\Bigr)\Bigr\|_{L^{1/(d-1)}(B(0,1))}\leq C\delta^{-\epsilon}\prod_{j=1}^{d}(\delta^{d-1}\: \mbox{\#}\tubes_{j}).
\end{equation*}
\end{theorem}

\vspace{3mm}

\begin{theorem}[Near-optimal multilinear restriction]\label{restthm}
For each $\epsilon> 0$, $q\geq\tfrac{2d}{d-1}$ and
$p'\leq \tfrac{d-1}{d}q$, there exists a constant $C$,
depending only on $A$, $\nu$, $\epsilon$,
$d$, $p$ and $q$, for which
$$
\Bigl\|\prod_{j=1}^{d}\ex_{j}g_{j}\Bigr\|_{L^{q/d}(B(0,R))}\leq C
R^{\epsilon}\prod_{j=1}^{d} \|g_{j}\|_{L^p(U_j)}
$$
for all $g_{j}\in L^p(U_j)$, $1\leq j\leq d$,
and all $R\geq 1$.
\end{theorem}

\vspace{3mm}

In Section 2 we show that Theorem \ref{restthm} follows from
Theorem \ref{kakthm}. In Section 4 we prove Theorem \ref{kakthm},
and in Section 5 we discuss the applications of our techniques to
lower orders of multilinearity, and to more general multilinear
$k$-plane transforms. In Section 6 we derive the natural variable
coefficient extensions of our results using further bootstrapping
arguments closely related to those of Bourgain. Finally, in Section
7 we give an application of our results to a variant of the
classical ``joints" problem considered in \cite{chaz}, \cite{sharir} and
\cite{feldman}.

\begin{remark}
The monotonicity approach that we take here arose from an attempt
to devise a continuous and more efficient version of an existing
induction-on-scales argument\footnote{This induction-on-scales
argument is closely related to that of Bourgain, and plays an
important role in our applications in Section 6.} introduced by
Wolff (and independently by the third author). This inductive
argument allows one to deduce linear (and multilinear) Kakeya
estimates for families of $\delta$-tubes from corresponding ones
for families of $\sqrt{\delta}$-tubes. However, unfortunately
there are inefficiencies present which prevent one from keeping
the constants in the inequalities under control from one iteration
to the next. Our desire to minimise these inefficiencies lead to
the introduction of the formulations in terms of gaussians adapted
to tubes (rather than rough characteristic functions of tubes).
The suggestion that one might then proceed by an
induction-on-scales argument, incurring constant factors of at
most $1$ at each scale, is then tantamount to a certain
monotonicity property. We should emphasise, however, that this
reasoning served mainly as philosophical motivation, and that
there are important differences between the arguments presented
here and the aforementioned induction arguments.
(Curiously however, one of
the most natural seeming formulations of monotonicity
fails at the endpoint $q=d/(d-1)$ when $d\geq 3$ -- see Proposition
\ref{notmon}.)  We also remark that closely related monotonicity arguments for
gaussians are effective in analysing
the Brascamp-Lieb inequalities \cite{CLL}, \cite{BCCT}; see
Remark \ref{brascamp-remark} below.
\end{remark}
\begin{remark}
There is perhaps some hope that variants of the techniques that we
introduce here may lead to progress on the original linear form of
the Kakeya conjecture. For $d\geq 3$ it seems unlikely that our
multilinear estimates (in their current forms) may simply be
``reassembled" in order to achieve this. However, our multilinear
results do suggest (in some non-rigorous sense) that if there were
some counterexamples to either the linear restriction or Kakeya
conjectures, then they would have to be somewhat ``non-transverse"
(or ``plany", in the terminology of \cite{KLT}). Issues of this nature
arise in our application to ``joints'' problems in Section \ref{joints-sec}.
\end{remark}
\begin{remark}
The $d$-linear transversality condition \eqref{spanning} has
also turned out
to be decisive when estimating
spherical averages
of certain multilinear extension operators
$(g_1,\hdots,g_d)\mapsto\ex_1 g_1\cdots\ex_d g_d$ of the type considered
here. See \cite{BBC} for further details.
\end{remark}
\subsubsection*{Notation}
For non-negative quantities $X$ and $Y$, we will use the statement
$X\lesssim Y$ to denote the existence of a constant $C$ for which
$X\leq CY$. The dependence of this constant on various parameters
will depend on the context, and will be clarified where
appropriate.

\subsection*{Acknowledgement}
We would like to thank Jim Wright for many helpful discussions on
a variety of techniques touched on in this paper.

\section{Multilinear Restriction $\iff$ Multilinear Kakeya}\label{RK}

It will be convenient to introduce some notation.
For $\alpha\geq 0$, $q\geq\tfrac{2d}{d-1}$ and
$p'\leq \tfrac{d-1}{d}q$, we use
$$\re^*(p\times\cdot\cdot\cdot\times p\rightarrow q;\alpha)$$
to denote the estimate
$$
\Bigl\|\prod_{j=1}^{d}\ex_{j}g_{j}\Bigr\|_{L^{q/d}(B(0,R))}\leq C
R^{\alpha}\prod_{j=1}^{d} \|g_{j}\|_{L^p(U_j)}
$$
for some constant $C$, depending only on $A$, $\nu$, $\alpha$,
$d$, $p$ and $q$, for all $g_{j}\in L^p(U_j)$, $1\leq j\leq d$,
and all $R\geq 1$. Similarly, for $\tfrac{d}{d-1}\leq
q\leq\infty$, we use
$$\ka^{*}(1\times\cdot\cdot\cdot\times 1\rightarrow q;\alpha)$$
to denote the estimate
\begin{equation}\label{Kalpha}
\Bigl\|\prod_{j=1}^{d}\Bigl(\sum_{T_{j}\in\tubes_{j}}
\chi_{T_{j}}\Bigr)\Bigr\|_{q/d}\leq C \delta^{-\alpha}
\prod_{j=1}^{d}(\delta^{d/q}\: \mbox{\#}\tubes_{j})
\end{equation}
for some constant $C$, depending only on $\alpha$, $d$ and $q$,
for all transversal collections of families of
$\delta$-tubes in
$\R^{d}$, and all $0<\delta\leq 1$. We again note that
\eqref{Kalpha} is equivalent by standard density arguments (in suitable weak topologies)
to the superficially stronger
\begin{equation}\label{Kalphasuper}
\Bigl\|\prod_{j=1}^{d}\Bigl(\sum_{T_{j}\in\tubes_{j}}\chi_{T_{j}}*\mu_{T_{j}}\Bigr)\Bigr\|_{L^{q/d}({\bf
R}^d)} \leq C \delta^{-\alpha}\prod_{j=1}^{d}\Bigl(\delta^{d/q}\:
\sum_{T_{j}\in\tubes_{j}}\|\mu_{T_{j}}\|\Bigr)
\end{equation}
for all finite measures $\mu_{T_{j}}$ ($T_j\in\tubes_j$,
$1\leq j\leq d$) on $\R^d$.

With this notation, Theorem \ref{kakthm} is equivalent to the statements
$\ka^{*}(1\times\cdot\cdot\cdot\times 1\rightarrow
q;0)$ for all $\tfrac{d}{d-1}<q\leq\infty$, and $\ka^{*}(1\times\cdot\cdot\cdot\times 1\rightarrow
\tfrac{d}{d-1};\epsilon)$ for all $\epsilon>0$;
Theorem \ref{restthm} is equivalent to
$\re^*(2\times\cdot\cdot\cdot\times 2\rightarrow
\tfrac{2d}{d-1};\epsilon)$ for all $\epsilon>0$.

As we have already discussed, a standard Rademacher-function
argument allows one to deduce the multilinear Kakeya conjecture
from the multilinear restriction conjecture; the linear argument found for instance in
\cite{Wolffbook} extends to the multilinear setting in a completely routine manner
and will not be detailed here. In the localised
setting this of course continues to be true; i.e. for any
$\alpha\geq 0$,
\begin{equation}\label{equ}
\re^*(2\times\cdot\cdot\cdot\times 2\rightarrow
\tfrac{2d}{d-1};\alpha) \implies
\ka^{*}(1\times\cdot\cdot\cdot\times 1\rightarrow
\tfrac{d}{d-1};2\alpha).
\end{equation}

Multilinearising a well-known bootstrapping argument of Bourgain
\cite{Bo} (again see \cite{TVV} for this argument in the bilinear
setting) we shall obtain the following reverse mechanism.
\begin{proposition}\label{MRiffMK}
For all $\alpha,\epsilon\geq 0$ and $\tfrac{2d}{d-1}\leq q\leq\infty$,
$$\re^*(2\times\cdot\cdot\cdot\times 2\rightarrow q;\alpha)\;+\;
\ka^{*}(1\times\cdot\cdot\cdot\times 1\rightarrow
\tfrac{q}{2};\epsilon)
\implies
\re^*(2\times\cdot\cdot\cdot\times 2\rightarrow
q;\tfrac{\alpha}{2}+\tfrac{\epsilon}{4}).
$$
\end{proposition}
Using elementary estimates we may
easily verify $\re^*(2\times\cdot\cdot\cdot\times 2\rightarrow
\tfrac{2d}{d-1};\alpha)$ for some large positive value of
$\alpha$. For example, noting that $|B(0,R)|=c_{d}R^{d}$ for some
constant $c_{d}$, we have that
$$
\Bigl\|\prod_{j=1}^{d}\ex_{j}g_{j}\Bigr\|_{L^{2/(d-1)}(B(0,R))}\leq
c_{d} R^{d(d-1)/2}\prod_{j=1}^{d} \|\ex_{j}g_{j}\|_{\infty} \leq
c_{d} R^{d(d-1)/2}\prod_{j=1}^{d}\|g_{j}\|_{L^{1}(U_{j})},
$$
which by the Cauchy--Schwarz inequality yields $\re^*(2\times\cdot\cdot\cdot\times 2\rightarrow
\tfrac{2d}{d-1};\tfrac{d(d-1)}{2})$.
In the presence of appropriately favourable Kakeya
estimates, this value of $\alpha$ may then be reduced by a
repeated application of the above proposition. In particular,
Proposition \ref{MRiffMK}, along with implication \eqref{equ},
easily allows one to deduce the equivalence
\begin{equation*}
\re^*(2\times\cdot\cdot\cdot\times 2\rightarrow
\tfrac{2d}{d-1};\epsilon)\forall\epsilon>0 \iff
\ka^{*}(1\times\cdot\cdot\cdot\times 1\rightarrow
\tfrac{d}{d-1};\epsilon)\forall\epsilon>0.
\end{equation*}
Arguing in much the same way allows us to reduce the proof of Theorem
\ref{restthm} to Theorem \ref{kakthm}, as claimed in the introduction.

The proof we give of Proposition \ref{MRiffMK} is very similar to
that of Lemma 4.4 of \cite{TVV}, and on a technical level is slightly
more straightforward. We begin by stating a lemma
which, given Remark \ref{alt} and the control of $|Y_{j}|$
implicit in \eqref{spanning} and \eqref{smoothcond}, is a standard
manifestation of the uncertainty principle (see \cite{Co} for the
origins of this, and
Proposition 4.3 of \cite{TVV} for a proof in the
bilinear case which immediately generalises to the multilinear case).
\begin{lemma}\label{uncertainty}
$\re^*(2\times\cdot\cdot\cdot\times 2\rightarrow
q;\alpha)$ is true if and only if
\begin{equation}\label{smudge}
\Bigl\|\prod_{j=1}^{d}\widehat{f}_{j}\Bigr\|_{L^{q/d}(B(0,R))}
\leq C R^{\alpha-d/2}\prod_{j=1}^{d}\|f_j\|_{2}
\end{equation}
for all
$R\geq 1$ and functions $f_j$ supported on
$A_{j}^{R}:=\Sigma_{j}(U_{j})+O(R^{-1})$, $1\leq j\leq d$.
\end{lemma}
We now turn to the proof of the proposition, where the implicit
constants in the $\lesssim$ notation will depend on at most
$A$, $\nu$, $d$, $p$, $\alpha$ and $\epsilon$. From the above
lemma it suffices to show that
$$
\Bigl\|\prod_{j=1}^{d}\widehat{f}_{j}\Bigl\|_{L^{q/d}(B(0,R))}
\lesssim R^{\alpha/2+\epsilon/4-d/2}\prod_{j=1}^{d}\|f_{j}\|_{L^2(A^R_j)}
$$
for all $f_j$ supported in $A_{j}^{R}$, $1\leq j\leq d$. To this
end we let $\phi$ be a real-valued bump function adapted to
$B(0,C)$, such that its Fourier transform is non-negative on the
unit ball. For each $R\geq 1$ and $x\in\R^d$ let
$\phi_{R^{1/2}}^{x}(\xi) := e^{-2\pi
ix\cdot\xi}R^{d/2}\phi(R^{1/2}\xi)$. Observe that
$\phi_{R^{1/2}}^{x}$ is an $L^1$-normalised modulated bump
function adapted to $B(0,C/R^{1/2})$, whose Fourier transform is
non-negative and bounded below on $B(x,R^{1/2})$, uniformly in
$x$. From the hypothesis $\re^*(2\times\cdot\cdot\cdot\times
2\rightarrow q;\alpha)$ and Lemma \ref{uncertainty}
we have
\begin{equation}\label{48}
\Bigl\|\prod_{j=1}^{d}\widehat{\phi_{R^{1/2}}^{x}}\widehat{f}_{j}\Bigr\|
_{L^{q/d}(B(x,R^{1/2}))}
\lesssim
R^{\alpha/2-d/4}\prod_{j=1}^{d}\|f_{j}*\phi_{R^{1/2}}^{x}\|_{L^2(A^R_j)}
\end{equation}
for all $x$. Averaging this over $x\in B(0,R)$ we obtain
\begin{equation}\label{49}
\Bigl\|\prod_{j=1}^{d}\widehat{f}_{j}\Bigl\|_{L^{q/d}(B(0,R))}
\lesssim R^{\alpha/2-d/4}\Bigl(
R^{-d/2}\int_{B(0,R)}\Bigl(\prod_{j=1}^{d}
\|f_{j}*\phi_{R^{1/2}}^{x}\|_{L^2({\bf R}^d)}^{2}\Bigr)^{q/(2d)}dx\Bigr)^{d/q}.
\end{equation}
Now for each $1\leq j\leq d$ we cover $A_{j}^{R}$ by a boundedly overlapping
collection of discs $\{\rho_{j}\}$ of diameter $R^{-1/2}$, and
set $f_{j,\rho_j}:=\chi_{\rho_{j}}f_{j}$. Since (for each $j$) the supports
of the functions $f_{j,\rho_j}*\phi_{R^{1/2}}^{x}$ have bounded overlap, it
suffices to show that
\begin{equation}\label{50}
\Bigl\|\prod_{j=1}^{d}\widehat{f}_{j}\Bigl\|_{L^{q/d}(B(0,R))}
\lesssim R^{\alpha/2-d/4}\Bigl(
R^{-d/2}\int_{B(0,R)}\Bigl(\prod_{j=1}^{d}
\sum_{\rho_{j}}\|f_{j,\rho_j}*\phi_{R^{1/2}}^{x}\|_{L^2({\bf R}^d)}^{2}
\Bigr)^{q/(2d)}
dx\Bigr)^{d/q}.
\end{equation}
The function $\widehat{\phi_{R^{1/2}}^{x}}$ is rapidly decreasing away from
$B(x,R^{1/2})$, and so by Plancherel's theorem, the left hand side
of \eqref{50} is bounded by
$$
C R^{\alpha/2-d/4}\Bigl(
R^{-d/2}\int_{B(0,R)}\Bigl(\prod_{j=1}^{d}
\sum_{\rho_{j}}\|\widehat{f_{j,\rho_j}}\|_{L^{2}(B(x,R^{1/2}))}^{2}
\Bigr)^{q/(2d)}dx\Bigr)^{d/q},
$$
since the portions of $\widehat{\phi_{R^{1/2}}^{x}}$ on translates of
$B(x,R^{1/2})$ can be handled by translation symmetry.
For each $\rho_j$ let $\psi_{\rho_j}$ be a Schwartz
function which is comparable
to $1$ on $\rho_j$ and whose Fourier transform satisfies
$$
|\widehat{\psi_{\rho_j}}(x+y)|\lesssim R^{-(d+1)/2}\chi_{\rho_{j}^*}(x)
$$
for all $x,y \in \R^d$ with $|y| \leq R^{1/2}$, where $\rho_{j}^*$ denotes an $O(R)\times
O(R^{1/2})\times \cdots\times O(R^{1/2})$-tube, centred at the
origin, and with long side pointing in the direction normal to the
disc $\rho_{j}$; the implicit constants here depending only on
$A$, $\nu$ and $d$. We point out that this is where we use the
full $C^{2}(U_j)$ control given by condition \eqref{smoothcond}.
If we define
$\tilde{f}_{j,\rho_j}:=f_{j,\rho_j}/\psi_{\rho_j}$, then $f_{j,\rho_j}$
and $\tilde{f}_{j,\rho_j}$ are pointwise comparable, and furthermore
by Jensen's inequality,
$$|\widehat{f}_{j,\rho_j}(x+y)|^{2}=|\widehat{\tilde{f}}_{j,\rho_{j}}*
\widehat{\psi}_{\rho_j}(x+y)|^{2}\lesssim
R^{-(d+1)/2}|\widehat{\tilde{f}}_{j,\rho_{j}}|^{2}*\chi_{\rho_{j}^*}(x)
$$
whenever $x \in \R^d$ and $|y| \leq R^{1/2}$.  Integrating this in $y$ we conclude
$$
\|\widehat{f}_{j,\rho_j}\|_{L^{2}(B(x,R^{1/2}))}^{2}\lesssim
R^{-1/2}|\widehat{\tilde{f}}_{j,\rho_{j}}|^{2}*\chi_{\rho_{j}^*}(x),
$$
and hence by rescaling the hypothesis $\ka(1\times\cdots\times 1\rightarrow
\tfrac{q}{2};\epsilon)$ (in its equivalent form \eqref{Kalphasuper}) we
obtain
\begin{eqnarray*}
\begin{aligned}
\Bigl\|\prod_{j=1}^{d}\widehat{f}_{j}&\Bigl\|_{L^{q/d}(B(0,R))}\\
&\lesssim
R^{\alpha/2-d/4}\Bigl(
R^{-d/2}\int_{B(0,R)}\Bigl(\prod_{j=1}^{d}
\sum_{\rho_{j}}R^{-1/2}|\widehat{\tilde{f}}_{j,\rho_{j}}|^{2}*\chi_{\rho_{j}^*}(x)
\Bigr)^{q/(2d)}dx\Bigr)^{d/q}\\
&\lesssim R^{\alpha/2+\epsilon/4-d/2}\prod_{j=1}^{d}
\Bigl(\sum_{\rho_{j}}\|\tilde{f}_{j,\rho_{j}}\|_{L^2(A^R_j)}^{2}\Bigr)^{1/2}\\
&\lesssim R^{\alpha/2+\epsilon/4-d/2}\prod_{j=1}^{d}\|f_{j}\|_{L^2(A^R_j)}.
\end{aligned}
\end{eqnarray*}
In the last two lines we have used Plancherel's theorem,
disjointness, and the pointwise comparability of
$\tilde{f}_{\rho_j}$ and $f_{j,\rho_j}$. This completes the proof of
Proposition \ref{MRiffMK}.

\section{The Loomis--Whitney case}

As we have already discussed, the ``Loomis--Whitney case" of
Conjecture \ref{MLKC} corresponds to the situation where each
family $\tubes_j$ consists of translates of a fixed tube with
direction $e_j$; the transverality hypothesis allows us to assume (after a linear
change of variables) that
$e_1,\ldots,e_d$ is the standard orthonormal basis of $\R^d$. The gaussian reformulation of this inequality (or
rather the equivalent inequality \eqref{1tubes}) alluded to in the
introduction is now
\begin{equation}\label{lwreformulation}
\Bigl\|\prod_{j=1}^{d}\Bigl(\sum_{v_{j}\in\vi_{j}}e^{-\pi\langle
A_{j}^{0}(\cdot-v_{j}),(\cdot-v_{j})\rangle}\Bigr)\Bigr\|_{L^{1/(d-1)}({\bf
R}^{d})} \leq C\prod_{j=1}^{d}\mbox{\#}\vi_{j},
\end{equation}
where for each $1\leq j\leq d$, $A_j^0$ is the orthogonal
projection to the $j^{\th}$ coordinate hyperplane $\{x\in\R^d:x_j =
0\}$, and $\vi_j$ is an arbitrary finite subset of $\R^d$. The
matrix $A_j^0$, which is just the diagonal matrix whose diagonal
entries are all $1$ except for the $j^{\th}$ entry, which is zero, we
refer to as the \textit{$j^{\th}$ Loomis--Whitney matrix}.

We will actually consider a rather more general $n$-linear setup
where there are $n$ distinct matrices $A_j$, which are not
necessarily commuting, and no relation between $n$ and $d$ is assumed.
We adopt the notation that $A \leq_{pd} B$
if $B-A$ is positive semi-definite, and $A
>_{pd} B$ if $B-A$ is positive definite.
The observations that if $A \leq_{pd} B$, then $D^* A D
\leq_{pd} D^* B D$ for any matrix $D$, and also if $A \geq_{pd} B
>_{pd} 0$ then $B^{-1} \geq_{pd} A^{-1} >_{pd} 0$ (as can be seen
by comparing the norms $\langle Ax, x \rangle^{1/2}$ and $\langle
Bx, x \rangle^{1/2}$ on $\R^d$ and then using duality) will be useful at
the end of the proof of the next proposition.

\begin{proposition}\label{lwmon} Let $d,n\geq 1$ and
$A_1, \ldots, A_n$ be positive semi-definite real symmetric
$d\times d$ matrices. Let
$\mu_1, \ldots, \mu_n$ be finite compactly supported
positive Borel measures on $\R^d$.  For $t
\geq 0$, $x \in \R^d$, and $1 \leq j \leq n$ let $f_j(t,x)$ denote
the non-negative quantity
$$ f_j(t,x) := \int_{{\bf R}^d} e^{-\pi \langle A_j(x-v_jt), (x-v_jt) \rangle}\ d\mu_j(v_j).$$
We interpret $t$ as the ``time'' variable.
Then if $p=(p_{1},\hdots,p_n)\in (0,\infty)^{n}$
is such that $p_{1}A_1+\cdots+p_{n}A_n$ is
non-singular\footnote{Note that this non-singularity is actually independent
of $p$, and is equivalent to the statement that
$\cap_{j=1}^{n}\ker A_j=\{0\}$.  See \cite{BCCT} for further analysis of the condition \eqref{a-jab}.} and
\begin{equation}\label{a-jab}
A_1, \ldots, A_n \leq_{pd} p_1A_1 + \cdots + p_{n}A_n,
\end{equation}
the quantity $$
Q_p(t) := \int_{{\bf R}^{d}}\prod_{j=1}^n
f_j(t,x)^{p_{j}}dx$$ is non-increasing in time.
\end{proposition}
\begin{corollary}\label{usemon1}
Under the conditions of Theorem \ref{lwmon},
$$\int_{{\bf R}^{d}}\prod_{j=1}^{n}f_{j}(1,x)^{p_{j}}dx
\leq\int_{{\bf R}^{d}}\prod_{j=1}^{n}f_{j}(0,x)^{p_{j}}dx
=\left(\frac{1}{\det
A_*}\right)^{\frac{1}{2}}\prod_{j=1}^{n}\|\mu_{j}\|^{p_{j}},
$$
where $A_*:=p_1A_1+\ldots+p_nA_n$.
\end{corollary}
\begin{remark}
The Loomis--Whitney case \eqref{lwreformulation} (with $C=1$) now follows from
Corollary \ref{usemon1} on setting $n=d$, $A_j = A_j^0$, $p_j =
1/(d-1)$ for all $j$, and the measures $\mu_j$ to be arbitrary
sums of Dirac masses.
\end{remark}
\begin{remark}\label{brascamp-remark}
Variants of Proposition \ref{lwmon} and Corollary \ref{usemon1}, by the current
authors and M. Christ \cite{BCCT}, have recently lead to new
proofs of the fundamental theorem of Lieb \cite{Lieb} concerning
the exhaustion by gaussians of the Brascamp--Lieb inequalities (of
which the Loomis--Whitney inequality is an important special case).
Although our proof of Proposition \ref{lwmon} is rather less
direct than the one given in \cite{BCCT} (which is closely related
to the heat-flow approach of \cite{CLL}), it does seem to lend
itself much better to the perturbed situation, as we will discover
in the next section.
\end{remark}
\begin{remark}
In the statement of Proposition \ref{lwmon},
the $v_j$'s can be thought of as the velocities with which the gaussians slide
to the origin.
\end{remark}

The proof we give of Proposition \ref{lwmon} is rather unusual.
We begin by considering the integer exponent case when $p=(p_1,\hdots,p_n)\in\N^n$. By
multiplying out the $(p_j)^{\th}$ powers in the expression for $Q_p(t)$, and using
Fubini's theorem, we may
obtain an explicit formula for the time derivative $Q_p'(t)$. With some careful
algebraic and combinatorial manipulation, we are then able to rewrite this
expression in a way that \emph{makes sense} for $p\not\in\N^n$, and is
manifestly non-positive whenever
$A_1, \ldots, A_n \leq_{pd} p_1A_1 + \cdots + p_{n}A_n$. Finally,
we appeal to an extrapolation lemma (see the appendix) to
conclude that the formula must in fact also hold for $p\in(0,\infty)^{n}$.
It may also be interesting to consider this approach in the light of the notion
of a fractional cartesian product of a set (see \cite{Blei}).

\subsection*{Proof of Proposition \ref{lwmon}}
We begin by considering the case when $p\in\N^{n}$. Then the quantity $Q_p(t)$ defined in this
proposition can be expanded as
$$ \int_{{\bf R}^d} \int_{({\bf R}^d)^{p_{1}}}
\cdots\int_{({\bf R}^d)^{p_{n}}} e^{-\pi \sum_{j=1}^n
\sum_{k=1}^{p_{j}} \langle A_j(x-v_{j,k} t), (x-v_{j,k} t)
\rangle}\ \prod_{j=1}^n \prod_{k=1}^{p_{j}} d\mu_j(v_{j,k}) dx.$$
On completing the square we find that
$$ \sum_{j=1}^n \sum_{k=1}^{p_{j}}
\langle A_j(x-v_{j,k} t), (x-v_{j,k} t) \rangle =  \langle A_*(x -
\overline{v} t), (x-\overline{v} t) \rangle + \delta t^2,$$ where
$A_* := \sum_{j=1}^n p_{j}A_j$ is a positive definite matrix,
$\overline{v} := A_*^{-1} \sum_{j=1}^n A_j \sum_{k=1}^{p_{j}}
v_{j,k}$ is the weighted average velocity, and $\delta$ is the weighted variance of the velocity,
\begin{equation}\label{delt}
 \delta := \sum_{j=1}^n \sum_{k=1}^{p_{j}} \langle A_j v_{j,k}, v_{j,k} \rangle - \langle A_* \overline{v}, \overline{v} \rangle.
\end{equation}
Using translation invariance in $x$, we thus have
$$
 Q'_p(t) = -2\pi t \int_{{\bf R}^d} \int_{({\bf R}^d)^{p_{1}}}
\cdots\int_{({\bf R}^d)^{p_{n}}} \delta \prod_{j=1}^n
\prod_{k=1}^{p_{j}} e^{-\pi \langle A_j(x-v_{j,k} t), (x-v_{j,k}
t) \rangle}\ d\mu_j(v_{j,k})dx.$$ If for each $j$ we let $\v_j$ be $v_j$
regarded as a
random variable associated to the probability measure
$$ \frac{e^{-\pi \langle A_j(x-vt), (x-vt) \rangle}d\mu_j(v_j)}{f_j(t,x)},$$
and let $\v_{j,1}, \ldots, \v_{j,p_{j}}$ be $p_j$ independent
samples of these random variables (with the $\v_{j,k}$ being
independent in both $j$ and $k$), then we can write the above as
$$
 Q'_p(t) = -2\pi t \int_{{\bf R}^d}
\E(\delta) \prod_{j=1}^n f_{j}(t,x)^{p_{j}}\ dx
$$
where $\delta$ is now considered a function of the $\v_{j,k}$, and $\E()$ denotes probabilistic expectation.
By linearity of expectation we have
\begin{equation}\label{Edelt}
 \E(\delta) = \sum_{j=1}^n \sum_{k=1}^{p_{j}} \E( \langle A_j \v_{j,k}, \v_{j,k} \rangle )
- \E( \langle A_* \overline{\v}, \overline{\v} \rangle ),
\end{equation}
(where of course $\overline{\v}$ is $\overline{v}$ regarded as a random
variable).
By symmetry, the first term on the right-hand side is
$\sum_{j=1}^n p_{j} \E( \langle A_j \v_j, \v_j \rangle )$. As for
the second term, by definition of $\overline{\v}$ we have
\begin{align*}
 \E( \langle A_* \overline{\v}, \overline{\v} \rangle )
&= \E( \langle A_*^{-1} \sum_{j=1}^n \sum_{k=1}^{p_{j}} A_j \v_{j,k}, \sum_{j'=1}^n \sum_{k'=1}^{p_{j'}} A_{j'} \v_{j,k'} \rangle ) \\
&= \sum_{j=1}^n \sum_{j'=1}^n \sum_{k=1}^{p_{j}}
\sum_{k'=1}^{p_{j'}} \E( \langle A_*^{-1} A_j \v_{j,k}, A_{j'}
\v_{j,k'} \rangle ).
\end{align*}
When $(j,k) \neq (j',k')$ we can factorise the expectation
using independence and symmetry to obtain
\begin{align*}
 \E( \langle A_* \overline{\v}, \overline{\v} \rangle )
&= \sum_{j=1}^n p_{j} \E( \langle A_*^{-1} A_j \v_j, A_j \v_j  \rangle) \\
&\quad - \sum_{j=1}^n p_{j} \langle A_*^{-1} A_j \E(\v_j), A_j \E( \v_j ) \rangle \\
&\quad +  \sum_{1 \leq j, j' \leq n}p_{j}p_{j'} \langle A_*^{-1}
A_j \E(\v_j), A_{j'} \E(\v_{j'}) \rangle.
\end{align*}
Combining these observations together, we obtain
\begin{equation}\label{qp-form}
Q'_p(t) = -2\pi t \int_{{\bf R}^d} G(p, t, x) \prod_{j=1}^n
f_{j}(t,x)^{p_{j}}\ dx
\end{equation}
where $G$ is the function
\begin{align*}
G(p,t,x) &:= \sum_{j=1}^{n}
p_{j}\left\{\E(\langle(A_{j}-A_{j}A_{*}^{-1}A_{j})\v_{j},\v_{j}\rangle)
-\langle(A_{j}-A_{j}A_{*}^{-1}A_{j})\E(\v_{j}),\E(\v_{j})\rangle\right\}\\
&+ \sum_{j=1}^{n}p_{j}\langle
A_{j}\E(\v_{j}),\E(\v_{j})\rangle-\sum_{j,j'}
p_{j}p_{j'}\langle A_{*}^{-1}A_{j}\E(\v_{j}),A_{j'}\E(\v_{j'})\rangle\\
&=\sum_{j=1}^{n}
p_{j}\E(\langle(A_{j}-A_{j}A_{*}^{-1}A_{j})(\v_{j}-\E(\v_{j})),
(\v_{j}-\E(\v_{j}))\rangle\\
&+ \sum_{j=1}^{n}p_{j}\langle A_{j}(\E(\v_{j})-\E(\overline{\v})),
(\E(\v_{j})-\E(\overline{\v}))\rangle.
\end{align*}
Note that $G$ makes sense now not just for $p\in\N^n$, but for all
$p\in (0,\infty)^{n}$. On the other hand by the chain rule we have
\begin{equation}\label{train}
Q_{p}'(t)=2\pi\int_{{\bf R}^d}\Bigl(\sum_{k=1}^{n}
p_{k}\E\bigl(\langle A_{k}\v_{k},x-t\v_k\rangle\bigr)\Bigr)\prod_{j=1}^n f_{j}(t,x)^{p_{j}}\
dx.
\end{equation} Now we observe that
since the adjugate matrix $\adj(A_{*}):= \det(A_{*})A_*^{-1}$ is
polynomial in $p$, $\det(A_*)G(p,t,x)$ is also polynomial in $p$.
Hence multiplying both \eqref{qp-form} and \eqref{train} by
$\det(A_*)$ and using Lemma \ref{densitycor} of the appendix,
(along with the
hypothesis that $A_{*}$ is non-singular), we may deduce that
\eqref{qp-form} in fact holds for all $p\in (0,\infty)^{n}$. Now,
if we choose $p$ so that $A_*\geq_{pd} A_j$ holds for all $1\leq j\leq
n$, the inner product $\langle (A_j - A_j A_*^{-1} A_j) \cdot,
\cdot \rangle$ is positive semi-definite, and hence $G(p,t,x)$ is
manifestly non-negative for all $t,x$.  This proves
Proposition \ref{lwmon}.
\endprf

\begin{remark}
The second formula for $G(p,t,x)$ above may be re-expressed as
\begin{eqnarray*}
\begin{aligned}
G(p,t,x) &= \sum_{j=1}^{n}
p_{j}\E(\langle(A_{j}-A_{j}A_{*}^{-1}A_{j})(\v_{j}-\E(\v_{j})),
(\v_{j}-\E(\v_{j}))\rangle\\
\\&+ \sum_{j=1}^{n}p_{j}\langle A_{j}\E(\v_{j}),
\E(\v_{j})\rangle - \langle A_* \E(\overline{\v}),\E(\overline{\v})\rangle.
\end{aligned}
\end{eqnarray*}
While in this formulation it is not obvious at a glance that $G(p,t,x) \geq 0$,
we can make a ``centre of mass'' change in the preceeding argument to deduce
this. Indeed, if we subtract $v_0 = v_0(t,x)$ from each $v_{j,k}$ in the
definition \eqref{delt} of $\delta$, then the value of $\delta$ remains
unchanged. If we now choose $v_0 = \E(\overline{\v})$, the last
term in our expression for $G(p,t,x)$ vanishes, whence $G$ is
nonnegative as before.  This type of ``Galilean invariance'' will also be used crucially
in the next section to obtain a similar positivity.

While this proof of Proposition \ref{lwmon} is a little more laboured than the one we have presented
in \cite{BCCT} (see also \cite{CLL}), it nevertheless paves the way for the
arguments of the next section where the (constant) matrices $A_j$ are replaced
by \emph{random} matrices $\A_j$ and are thus subject to the expectation
operator $\E$. In that context, we are able to prove a suitable variant of the
formula for $G$ presented in this remark.
\end{remark}

\begin{remark} If $\sum_{j=1}^n p_j A_j > A_l$ for all $l$, we can expect
there to be room for a stronger estimate to hold. This will also play an
important role in the next section where we will use it to handle error terms
arising in our analysis.
\end{remark}

\section{The perturbed Loomis-Whitney case}

In this section we prove a perturbed version of Proposition
\ref{lwmon} of the previous section. Theorem \ref{kakthm} will
then follow as a special case. Although Theorem \ref{kakthm} is
our main goal, working at this increased level of generality has
the advantage of providing more general (multilinear) $k$-plane
transform estimates at all levels of multilinearity. We shall
discuss these further applications briefly in the next section.

If $A$ is a real symmetric $d \times d$ matrix, we use $\|A\|$ to denote the operator norm
of $A$ (one could also use other norms here, such as the Hilbert-Schmidt norm, as we are not tracking
the dependence of constants on $d$).

\begin{proposition}\label{wowgen}
Let $d,n\geq 1$, $\eps > 0$ and $M_{1},\hdots,M_{n}$ be positive
semi-definite real symmetric $d\times d$ matrices. In addition
suppose that
$p=(p_{1},\hdots,p_n)\in (0,\infty)^{n}$ is such that the sum
$p_{1}M_{1}+\cdots+p_{n}M_{n}$ is non-singular and
\begin{equation}\label{gap}
p_{1}M_{1}+\cdots+p_{n}M_{n}>_{pd} M_{j}
\end{equation}
for all $1\leq j\leq n$. For each $1\leq j\leq n$ let $\Omega_j$
be a collection of pairs $(A_j,v_j)$, where $v_j \in \R^d$, $A_j
\geq_{pd} 0$, $\|A_j^{1/2} - M_j^{1/2} \| \leq \eps$, and let
$\mu_j$ be a finite compactly supported positive Borel
measure on $\Omega_j$.  For $t \geq 0$, $x
\in \R^d$, and $1 \leq j \leq n$ let $f_j(t,x)$ denote the
non-negative quantity
$$ f_j(t,x) := \int_{\Omega_j} e^{-\pi \langle A_j(x-v_jt), (x-v_jt) \rangle}\ d\mu_j(A_j,v_j).$$
Then if $\eps$ is sufficiently small depending on $p$ and the
$M_{j}$'s, we have the approximate monotonicity formula
$$\int_{{\bf R}^{d}}\prod_{j=1}^n f_j(t,x)^{p_{j}}dx\leq (1 + O(\eps))
\int_{{\bf R}^d} \prod_{j=1}^n f_j(0,x)^{p_{j}}dx$$ for all $t
\geq 0$, where we use $O(X)$ to denote a quantity bounded by $C X$
for some constant $C> 0$ depending only on $p$ and the $M_{j}$'s.
\end{proposition}
\begin{corollary}\label{corwow}
If $p$ is such that \eqref{gap} holds, then if $\epsilon$ is small
enough depending on $p$ and the $M_j$'s, we have the inequality
$$\int_{{\bf R}^{d}}\prod_{j=1}^{n}f_{j}(1,x)^{p_{j}}dx
\leq (1+O(\epsilon))\left(\frac{1}{\det
M_*}\right)^{\frac{1}{2}}\prod_{j=1}^{n}\|\mu_{j}\|^{p_{j}},
$$
where $M_*=p_1M_1+\cdots+p_nM_n$.
\end{corollary}
To prove the corollary it is enough to show that
there is a $c>0$ (depending only on $p$ and the $M_j$'s) such that
\begin{equation}\label{labour}
\prod_{j=1}^{n}f_{j}(0,x)^{p_{j}}\leq
e^{-\pi\langle(1-c\epsilon)M_*x,x\rangle}\prod_{j=1}^{n}\|\mu_{j}\|^{p_{j}}
\end{equation}
for all $x\in\R^{d}$.
To this end we first observe that for any $\lambda>0$ we have the identity
\begin{eqnarray*}
\begin{aligned}
\prod_{j=1}^{n}f_{j}&(0,x)^{p_{j}}
=\prod_{j=1}^{n}\left(\int_{\Omega_{j}}e^{-\pi\langle A_jx,x
\rangle}d\mu_j(A_j,v_j)\right)^{p_j}\\
&=e^{-\pi\langle (1-\lambda\epsilon\sum_{k=1}^{n}p_k)M_*x,x\rangle}
\prod_{j=1}^{n}\left(\int_{\Omega_j}e^{-\pi\langle (A_j-M_j+\lambda\epsilon M_*)x,x
\rangle}d\mu_j(A_j,v_j)\right)^{p_j}.\\
\end{aligned}
\end{eqnarray*}
Now since $M_*>_{pd}0$, there exists a constant $c'>0$, such that
$M_{*}\geq_{pd}c'I$, where $I$ denotes the identity matrix. Using this and
\eqref{ajm-approx} below, we may choose $\lambda>0$ (depending only on
$p$ and the $M_j$'s) such that $A_j-M_j+\lambda\epsilon M_{*}$ is positive
definite, and hence $e^{-\pi\langle (A_j-M_j+\lambda\epsilon M_*)x,x\rangle}
\leq 1$ for all $x$. Setting $c=\lambda\sum p_k$ completes the proof of
\eqref{labour}, and hence the corollary.
\begin{remark}  When the directions of the tubes in $\tubes_j$
are sufficiently close to $e_j$, Theorem \ref{kakthm} follows from
Corollary \ref{corwow}
on setting $n=d$, $M_j=A_j^0$, $p_{j}=p>\tfrac{1}{d-1}$ and the
measures $\mu_{j}$ to be appropriate sums of Dirac masses.
By affine invariance and the triangle inequality one can
handle any linearly independent direction sets, although of
course the constants may now get significantly larger than 1.
\end{remark}
\subsection*{Proof of Proposition \ref{wowgen}} From the estimate $\|
A_j^{1/2} - M_j^{1/2} \| \leq \eps$
we have
\begin{equation}\label{ajm-sqrt}
 A_j^{1/2} = M_j^{1/2} + O(\eps)
\end{equation}
and thus
\begin{equation}\label{ajm-approx} A_j = M_j + O(\eps).
\end{equation}

Again, we begin by considering the case when $p$ is an integer.
Let $Q_p$ denote the quantity
$$ Q_p(t) := \int_{{\bf R}^{d}}\prod_{j=1}^n f_j(t,x)^{p_{j}}dx;$$ we can expand this as $Q_{p}(t)=$
$$\int_{{\bf R}^d} \int_{\Omega_{1}^{p_{1}}} \ldots \int_{\Omega_{n}^{p_{n}}}
e^{-\pi \sum_{j=1}^n \sum_{k=1}^{p_j} \langle A_{j,k}(x-v_{j,k}
t), (x-v_{j,k} t) \rangle}\ \prod_{j=1}^n \prod_{k=1}^{p_j}
d\mu_j(A_{j,k}, v_{j,k}) dx.$$ It turns out that this quantity
will not be easy for us to study directly (mainly because of the
quantity $A_*^{-1}$ which will appear in the derivative of $Q_p$).
Instead, we shall consider the modified quantity $\tilde Q_p(t)$
defined by $\tilde Q_p(t) :=$
\begin{equation}\label{qptm-def}
\int_{{\bf R}^d} \int_{\Omega_1^{p_1}} \cdot\cdot
\int_{\Omega_n^{p_n}} \det(A_*) e^{-\pi \sum_{j=1}^n
\sum_{k=1}^{p_j} \langle A_{j,k}(x-v_{j,k} t), (x-v_{j,k} t)
\rangle}\ \prod_{j=1}^n \prod_{k=1}^{p_j} d\mu_j(A_{j,k}, v_{j,k})
dx
\end{equation}
where $A_* := \sum_{j=1}^n \sum_{k=1}^{p_j} A_{j,k}$ is a positive
definite matrix: the point is that
the determinant $\det(A_*)$
will eventually be used to convert $A_*^{-1}$ into a quantity
which is a \emph{polynomial} in the $A_{j,k}$ (the adjugate or cofactor
matrix of $A_*$); i.e. $\det(A_*)A_*^{-1}=\adj(A_*)$.

Let us now see why the weight $\det(A_*)$ is mostly harmless. From
\eqref{ajm-approx} we have $A_* = M_* + O(\eps)$, where
$M_*=p_1M_{1}+\hdots+p_nM_{n}$, and so in particular we have
$$ \det(A_*) = \det(M_*) + O(\eps).$$
Thus for $\eps$ small enough, $\det(A_*)$ is comparable to the
positive constant $\det(M_*)$.  Since all the
terms in $Q_p(t)$ and $\tilde Q_p(t)$ are non-negative, we have
thus established the bound
\begin{equation}\label{qpm-qtpm}
 Q_p(t) = (1 + O(\eps)) \det(M_*)^{-1} \tilde Q_p(t)
\end{equation}
for \emph{integer} $p_1,\hdots,p_n$.  However, for applications we
need this relation for non-integer $p_1,\hdots p_n$.  There is
an obvious difficulty in doing so, namely that $\tilde Q_p(t)$ is
not even defined for $p\not\in\N^n$.  However this can be fixed by
performing some manipulations (similar to ones considered
previously) to rewrite $\tilde Q_p(t)$ as an expression which
makes sense for arbitrary $p\in(0,\infty)^n$.

To this end, for each $j$ we let $(\A_j, \v_j)$ be random
variables (as before) associated to the probability measure
$$ \frac{e^{-\pi \langle A_j(x-v_{j}t), (x-v_{j}t)
\rangle}d\mu_{j}(A_{j},v_{j})}{f_j(t,x)}$$ and let $(\A_{j,1},
\v_{j,1}), \ldots, (\A_{j,p_j}, \v_{j,p_j})$ be $p_j$ independent
samples of these random variables (with the $(\A_{j,k}, \v_{j,k})$
being independent of $(\A_{j',k'}, \v_{j',k'})$ when $(j,k) \neq
(j',k')$).  Then we can rewrite \eqref{qptm-def} as
$$
 \tilde Q_p(t) = \int_{{\bf R}^d}
\E(\det(\A_*)) \prod_{j=1}^n f_j(t,x)^{p_j}\ dx,
$$
where $\A_* := \sum_{j=1}^n \sum_{k=1}^{p_j} \A_{j,k}$.  If we
write $\R_{j,k} := M_j - \A_{j,k}$ then the random variables
$\R_{j,k}$ are $O(\eps)$ (by \eqref{ajm-approx}).
Observe that for fixed $j$, all the $\R_{j,k}$ have the same
distribution as some fixed random variable $\R_j$, and we can
write
$$ \det(\A_*) = \det(M_*) + P( (\R_{j,k})_{1 \leq j \leq n;
1 \leq k \leq p_j} )$$ where $P$ is a polynomial in the
coefficients of the $\R_{j,k}$ which has no constant term (i.e.
$P(0) = 0$).  Thus we have
$$ \tilde Q_p(t) = \det(M_*) Q_p(t) +
\int_{{\bf R}^d} \E(P(\R_{j,k})_{1 \leq j \leq n; 1 \leq k \leq
p_j}) \prod_{j=1}^n f_j(t,x)^{p_j}\ dx.$$ Now by taking advantage
of independence and symmetry of the random variables $\R_{j,k}$,
we can write the expression $\E(P(\R_{j,k}))$ as a polynomial
combination of $p$ and of the (tensor-valued) moments
$\E(\R_j^{\otimes m})$ for some finite number $m=1,\ldots,M$ of
$m$ (with $M$ depending only on $n$, $d$ and of course the
$M_{j}$'s); thus we have
\begin{equation}\label{qptm-form}  \tilde Q_p(t) = \det(M_*) Q_p(t) +
\int_{{\bf R}^d} \tilde P(p, \E(\R_{j}^{\otimes m})_{1 \leq j \leq
n; 1 \leq m \leq M}) \prod_{j=1}^n f_j(t,x)^{p_j}\ dx
\end{equation}
for some polynomial $\tilde P$ depending only on $n$, $d$ and the
$M_{j}$'s. Since $P$ has no constant term it is easy to see that
$\tilde P(p, \cdot)$ also has no constant term, i.e. $\tilde
P(p,0) = 0$ (this can also be seen by considering the $\eps=0$
case). We remark that this polynomial can be computed explicitly
using the Lagrange interpolation formula, although we make no use
of this here.
The right-hand side of \eqref{qptm-form} makes sense for any
$p_1,\hdots,p_n>0$, not necessarily integers, and so we shall
adopt it as our \emph{definition} of $\tilde Q_p$ in general. Since $\R_j
= O(\eps)$ and $\tilde P$ has no constant term we observe that
$\tilde P(p, \E(\R_{j}^{\otimes m})_{1 \leq j \leq n; 1 \leq m
\leq M}) = O(\eps)$, whence we obtain \eqref{qpm-qtpm} for all
$p_1,\hdots,p_n
> 0$ (not just the integers), though of course the implicit
constants in the $O$ notation will certainly depend (polynomially)
on $p$.

In order to prove the proposition, it thus suffices by
\eqref{qpm-qtpm} to show that the quantity $\tilde Q_p(t)$ is
non-decreasing in time for sufficiently small $\eps$.
Again, we begin by working with
$p\in\N^n$, so that we may use \eqref{qptm-def}. Now we
differentiate $\tilde Q_p(t)$.  As before, we can complete the square and
write
$$ \sum_{j=1}^n \sum_{k=1}^{p_j}
\langle A_{j,k}(x-v_{j,k} t), (x-v_{j,k} t) \rangle = \langle
A_*(x - \overline{v} t), (x-\overline{v} t) \rangle + \delta t^2$$
where $\overline{v} := A_*^{-1} \sum_{j=1}^n \sum_{k=1}^{p_j}
A_{j,k} v_{j,k}$ is the weighted average velocity, and $\delta$ is the weighted variance
\begin{equation}\label{delta}
 \delta := \sum_{j=1}^n \sum_{k=1}^{p_j} \langle A_{j,k} v_{j,k}, v_{j,k} \rangle -
\langle A_* \overline{v}, \overline{v} \rangle.
\end{equation}
Arguing as in the previous section, we thus have
$\tilde Q'_p(t) =$
$$-2\pi t \int_{{\bf R}^d} \int_{\Omega_1^{p_1}} \ldots \int_{\Omega_n^{p_n}}
\det(A_*) \delta \prod_{j=1}^n \prod_{k=1}^{p_j} e^{-\pi \langle
A_{j,k}(x-v_{j,k} t), (x-v_{j,k} t) \rangle}\ d\mu_j(A_{j,k},
v_{j,k})\ dx.$$ Recalling the random variables
$\A_{j,k}$, $\v_{j,k}$, we thus have
\begin{equation}\label{tqpm-deriv}
 \tilde Q'_p(t) = -2\pi t \int_{{\bf R}^d} \E(\det(\A_*) \delta) \prod_{j=1}^n f(t,x)^{p_j}\ dx
\end{equation}
where $\delta$ is now considered a function of the $\A_{j,k}$ and
$\v_{j,k}$.

Let us rewrite $\delta$ slightly by inserting the definition of
$\overline{\v}$ and using the self-adjointness of
$A_*^{-1}$, to obtain
$$ \delta = \sum_{j=1}^n
\sum_{k=1}^{p_j} \langle \A_{j,k} \v_{j,k}, \v_{j,k} \rangle -
\langle \A_*^{-1}\sum_{j=1}^n \sum_{k=1}^{p_j} \A_{j,k} \v_{j,k},
\sum_{j'=1}^n \sum_{k'=1}^{p_{j'}} \A_{j',k'} \v_{j',k'}\rangle.$$
We now take advantage of a certain ``Galilean invariance'' of the
problem.  We introduce an arbitrary (deterministic) vector field
$v_0(t,x)$ which we are at liberty to select later, and observe
that the above expression is unchanged if we replace all the
$\v_{j,k}$ by $\v_{j,k}-v_0$: \begin{eqnarray*}
\begin{aligned}\delta = \sum_{j=1}^n &\sum_{k=1}^{p_j} \langle \A_{j,k}
(\v_{j,k}-v_0), (\v_{j,k}-v_0) \rangle \\&- \langle \A_*^{-1}
\sum_{j=1}^n \sum_{k=1}^{p_j} \A_{j,k} (\v_{j,k}-v_0),
\sum_{j'=1}^n \sum_{k'=1}^{p_{j'}} \A_{j',k'}
(\v_{j',k'}-v_0)\rangle.
\end{aligned}
\end{eqnarray*}
If we multiply this quantity by $\det(\A_*)$ then we
obtain a polynomial:
\begin{eqnarray}\label{deltam-poly}
\begin{aligned}
 \det(\A_*) \delta = \sum_{j=1}^n &\sum_{k=1}^{p_j} \det(\A_*) \langle \A_{j,k} (\v_{j,k}-v_0), (\v_{j,k}-v_0) \rangle\\ &-
\langle \adj(\A_*)\sum_{j=1}^n \sum_{k=1}^{p_j} \A_{j,k}
(\v_{j,k}-v_0),\sum_{j'=1}^n \sum_{k'=1}^{p_{j'}} \A_{j',k'}
(\v_{j',k'}-v_0)\rangle
\end{aligned}
\end{eqnarray}
where $\adj(\A_*)$ is the adjugate or cofactor matrix of $\A_*$,
which is a polynomial in the coefficients of $\A_*$ and is thus
polynomial in the coefficients of the $\A_{j,k}$. To make this
expression more elliptic, we introduce the matrices $\B_{j,k} :=
\A_{j,k}^{1/2}$ and the vectors $\w_{j,k} := \B_{j,k}
(\v_{j,k}-v_0)$, and write \eqref{deltam-poly} as
\begin{eqnarray}\label{deltam-poly-2}
\begin{aligned}
\det(\A_*) \delta = \sum_{j=1}^n &\sum_{k=1}^{p_j} \det(\A_*)
\|\w_{j,k}\|^2\\& - \langle \adj(\A_*)\sum_{j=1}^n
\sum_{k=1}^{p_j} \B_{j,k} \w_{j,k}, \sum_{j'=1}^n
\sum_{k'=1}^{p_{j'}} \B_{j',k'} \w_{j',k'}\rangle.
\end{aligned}
\end{eqnarray}
Note that $\A_* = \sum_{j=1}^n \sum_{k=1}^{p_j} \B_{j,k}^2$ is a
polynomial in the $\B_{j,k}$, so the right-hand side of
\eqref{deltam-poly-2} is a polynomial in the $\B_{j,k}$ and
$\w_{j,k}$. We can then take expectations, taking advantage of the
independence and symmetry of the random variables $(\B_{j,k},
\w_{j,k})$, and obtain a formula of the form
\begin{equation}\label{Sm-form}
 \E(\det(\A_*) \delta) = S(p, (\E( (\B_j, \w_j)^{\otimes m} ))_{1 \leq j \leq n; 1 \leq m \leq M} )
\end{equation}
for some polynomial $S$, and where the power of moments $M$
depends only on $d$ (in fact it is $2d+2$, which is the degree of
the polynomial in \eqref{deltam-poly-2}).  Here of course $(\B_j,
\w_j)$ represents any random variable with the same distribution
as the $(\B_{j,k}, \w_{j,k})$.  Also note that as the right-hand
side of \eqref{deltam-poly-2} is purely quadratic in the
$\w_{j,k}$, the expression \eqref{Sm-form} must be purely quadratic in the
$\w_j$. Just as before, one may of course use the Lagrange
interpolation formula here to write down an explicit expression
for $S$.

Inserting this formula back into \eqref{tqpm-deriv} we obtain
\begin{equation}\label{tqp-deriv-poly}
\tilde Q'_p(t) = -2\pi t \int_{{\bf R}^d}
S(p, (\E( (\B_j, \w_j)^{\otimes m} ))_{1 \leq j \leq n; 1 \leq m
\leq M} ) \prod_{j=1}^n f(t,x)^{p_j}\ dx
\end{equation}
when $p\in\N^n$.  However, since
\begin{equation*}
Q_{p}'(t)=2\pi\int_{{\bf R}^d}\Bigl(\sum_{k=1}^{n}
p_{k}\E\bigl(\langle
\A_{k}\v_{k},x-t\v_k\rangle\bigr)\Bigr)\prod_{j=1}^n f_{j}(t,x)^{p_{j}}\
dx,
\end{equation*}
by Lemma \ref{densitycor} of the appendix,
identity \eqref{tqp-deriv-poly}
must also be true for arbitrary $p\in (0,\infty)^n$.  Thus to
conclude the proof it will suffice to show that
$$ S(p, (\E( (\B_j, \w_j)^{\otimes m} ))_{1 \leq j \leq n;
1 \leq m \leq M} ) \geq 0$$
for any $p$ satisfying the hypotheses of the theorem, if $\eps$ is
sufficiently small depending on $p$ and the $M_j$'s.

From \eqref{ajm-sqrt}
we know that $\B_j = M_{j}^{1/2} + O(\eps)$, and in particular we have
$\B_j = O(1)$.  Since $S$ is purely quadratic in the $\w_j$, we
thus have
\begin{eqnarray}\label{tony}
\begin{aligned}
S(p, (\E( (\B_j, \w_j)^{\otimes m} ))_{1 \leq j \leq n; 1 \leq m
\leq M} ) &= S(p, (\E( (M_j^{1/2}, \w_j)^{\otimes m} ))_{1 \leq j
\leq n; 1 \leq m \leq M} )\\& - O\bigl(\eps \sum_{j=1}^n \E(
\|\w_j\|^2 )\bigr),
\end{aligned}
\end{eqnarray}
as we can use the Cauchy--Schwarz inequality to control any
cross terms such as $\E( \w_j )^{\otimes 2}$ or $\E(\w_j) \otimes
\E(\w_k)$.

\begin{lemma}\label{formulaforS}
For arbitrary $p_1,\hdots,p_n>0$
\begin{equation}\label{equationforS}
\begin{split}
S(p, (\E(& (M_j^{1/2}, \w_j)^{\otimes m} ))_{1 \leq j \leq n; 1
\leq m \leq M} )\\&=
\det M_* \Bigl\{\sum_{j=1}^n p_{j}\|\E(\w_j)\|^2\\
&+ \sum_{j=1}^n p_{j}\E (\langle(I-M_j^{1/2}M_*^{-1}M_j^{1/2})(\w_j-\E(\w_j)),
(\w_j-\E(\w_j))\rangle) \\
&-\| M_*^{-1/2}\sum_{j=1}^n p_{j}M_j^{1/2} \E(\w_j) \|^2\Bigr\}\\
\end{split}
\end{equation}
\end{lemma}

\begin{proof}
As usual we begin with integer $p_1, \hdots, p_n$.
From \eqref{Sm-form} and
\eqref{deltam-poly-2} (with the $\B_j$ replaced by $M_j^{1/2}$) we
have
\begin{align*} S(p, (\E( (M_j^{1/2}, &\w_j)^{\otimes m} ))_{1 \leq j \leq n; 1 \leq m \leq M} )
= \E\bigl(\sum_{j=1}^n \sum_{k=1}^{p_j} \det(M_*) \|\w_{j,k}\|^2 \\
&- \langle \adj(M_*) \sum_{j=1}^n \sum_{k=1}^{p_j} M_{j}^{1/2}
\w_{j,k}, \sum_{j'=1}^n \sum_{k'=1}^{p_{j'}} M_{j'}^{1/2}
\w_{j',k'} \rangle\bigr).
\end{align*}
Writing $\adj(M_*) = \det(M_*) (M_*)^{-1}$, we thus have
\begin{align*}
 S(p, (\E( (M_j^{1/2}, \w_j)^{\otimes m} ))&_{1 \leq j \leq n; 1 \leq m \leq M} )\\& =
 \det(M_*)\sum_{j=1}^n \sum_{k=1}^{p_j} \E( \| \w_{j,k} \|^2 )\\
&- \det(M_*)\E (\langle M_{*}^{-1}\sum_{j=1}^n \sum_{k=1}^{p_j}
M_{j}^{1/2} \w_{j,k}, \sum_{j'=1}^n \sum_{k'=1}^{p_{j'}}
M_{j'}^{1/2} \w_{j',k'} \rangle ),
\end{align*}
which by independence and symmetry we can rewrite as
\begin{eqnarray*}
\begin{aligned}
\det(M_*)&\sum_{j=1}^n p_{j}\E( \| \w_j \|^2 )\\
&-\det(M_*)\sum_{j=1}^{n} p_{j}\E( \langle M_{*}^{-1}
M^{1/2}_j \w_j, M^{1/2}_j \w_j \rangle) \\
&- \det(M_*) \sum_{j=1}^{n} p_{j}(p_{j}-1)\langle
M_{*}^{-1}M_j^{1/2} \E(\w_j), M_j^{1/2}
\E(\w_j) \rangle\\
&- 2\det(M_*) \sum_{1 \leq j < j' \leq n} p_{j}p_{j'}\langle
M_{*}^{-1} M_j^{1/2} \E(\w_j),
M_{j'}^{1/2} \E(\w_{j'}) \rangle.
\end{aligned}
\end{eqnarray*}
Now we argue as in the previous section and write
\begin{eqnarray*}
\begin{aligned}
\sum_{j=1}^n p_{j}\E( \| \w_j \|^2 ) - \sum_{j=1}^n p_{j} &\E(
\langle M_*^{-1}M_j^{1/2} \w_j, M_j^{1/2} \w_j \rangle)\\
&= \sum_{j=1}^n p_{j}\E( \langle
(I-M_j^{1/2}M_{*}^{-1}M_{j}^{1/2})\w_j,\w_j\rangle ),
\end{aligned}
\end{eqnarray*}
and
\begin{eqnarray*}
\begin{aligned} \E( \langle (I-&M_j^{1/2}M_{*}^{-1}M_{j}^{1/2}) \w_j, \w_j \rangle )\\
&= \langle (I-M_j^{1/2}M_{*}^{-1}M_{j}^{1/2}) \E(\w_j), \E(\w_j) \rangle\\& + \E(
\langle (I-M_j^{1/2}M_{*}^{-1}M_{j}^{1/2}) (\w_j - \E(\w_j)), (\w_j - \E(\w_j))
\rangle ).
\end{aligned}
\end{eqnarray*}

Hence
\begin{eqnarray*}
\begin{aligned}
S(p, (\E( (M_j^{1/2}, \w_j&)^{\otimes m} ))_{1 \leq j
\leq n; 1 \leq m \leq M} )\\
= \det (M_*) &\Bigl\{ \sum_{j=1}^n p_j\langle
(I-M_j^{1/2}M_{*}^{-1}M_{j}^{1/2}) \E(\w_j), \E(\w_j)
\rangle\\
&+ \sum_{j=1}^n p_j \E(\langle (I-M_j^{1/2} M_* ^{-1}
M_j^{1/2})(\w_j - \E(\w_j)),
(\w_j - \E(\w_j)) \rangle\\
&-\sum_{j=1}^{n}p_{j}(p_{j}-1)\langle M_*^{-1}M_j^{1/2}\E(\w_j), M_j^{1/2} \E(\w_j) \rangle\\
&- 2\sum_{1 \leq j < j' \leq n} p_{j}p_{j'}\langle
M_*^{-1}M_j^{1/2}\E(\w_j),
M_{j'}^{1/2} \E(\w_{j'}) \rangle\Bigr\}.
\end{aligned}
\end{eqnarray*}

We can rearrange the right-hand side as
\begin{eqnarray*}
\begin{aligned}
\det (M_*) \Bigl\{&\sum_{j=1}^n p_{j}\|\E(\w_j)\|^2 \\&+
\sum_{j=1}^n p_j \E(\langle (I-M_j^{1/2} M_* ^{-1} M_j^{1/2})
(\w_j - \E(\w_j)),
(\w_j - \E(\w_j)) \rangle\\
&-\sum_{j=1}^{n}p_{j}^2\langle M_*^{-1}M_j^{1/2} \E(\w_j),
M_j^{1/2}\E(\w_j)
\rangle\\
&- 2\sum_{1 \leq j < j' \leq d} p_{j}p_{j'}\langle
M_*^{-1}M_j^{1/2}\E(\w_j),
M_{j'}^{1/2} \E(\w_{j'}) \rangle \Bigr\}
\end{aligned}
\end{eqnarray*}
which can be rearranged further as
\begin{align*}
\det (M_*) \Bigl\{&\sum_{j=1}^n p_{j}\|\E(\w_j)\|^2 \\&+
\sum_{j=1}^n p_j \E(\langle (I-M_j^{1/2} M_* ^{-1} M_j^{1/2})(\w_j
- \E(\w_j)),
(\w_j - \E(\w_j)) \rangle\\
&-\| M_*^{-1/2}\sum_{j=1}^n p_{j}M_j^{1/2} \E(\w_j) \|^2 \Bigr\}
\end{align*}
as desired when the $p_j$ are integers. Applying Lemma
\ref{densitycor} of the appendix (in the now familiar way)
concludes the proof of the lemma.
\end{proof}

The last term in \eqref{equationforS} has an unfavourable sign.
However, we can at last select our vector field $v_0$ in order to remove this term.
More precisely, we choose $v_0$ by the formula
$$v_{0} := \E\bigl(\sum p_jM_j^{1/2}\A_{j}^{1/2}\bigr)^{-1}
\sum p_{j}\E(M_{j}^{1/2}\A_{j}^{1/2}\v_{j}),$$
which is well defined for $\eps$ sufficiently small since $\A_j = M_j + O(\eps)$ and
thus $\sum p_jM_j^{1/2}\A_{j}^{1/2}$ is close to the invertible matrix $M_*$.  With this choice of $v_0$
we compute that
$$ \sum_{j=1}^n p_{j}M_j^{1/2} \E(\w_j) = 0$$
and so the last term in \eqref{equationforS} vanishes.

Recalling \eqref{tony}, \eqref{equationforS}, and using the fact that there exists a positive
constant $c$ depending only on $p$ and the $M_j$'s, such that
$I-M_j^{1/2}M_{*}^{-1}M_{j}^{1/2} \geq_{pd} c I$
for all $1\leq j\leq n$, we see that it suffices to show that
$$
\sum_{j=1}^n p_{j}\|\E( \w_j)\|^2  + c\sum_{j=1}^n p_{j}\E( \|
\w_j - \E(\w_j) \|^2 ) - O\bigl(\eps \sum_{j=1}^n \E( \|\w_j\|^2
)\bigr) \geq 0.
$$
However we have
$\| \w_j \|^2 \leq 2\| \w_j - \E(\w_j) \|^2 + 2\| \E(\w_j) \|^2,$
and so the claim now follows if $\eps$ is sufficiently small
depending on $p$ and the $M_j$'s.

\vspace{2mm}

When each $p_j$ is an integer, the quantity $Q_p (t)$ is monotonic
without further hypotheses on the matrices $A_j$ since $\delta$ defined
in \eqref{delta} is manifestly non-negative. This in particular includes
the two dimensional bilinear Kakeya situation $n=d=2$ and $p_1=p_2=1$.
(While we do not in this paper establish monotonicity of $Q_{p}(t)$ under
the hypotheses of Proposition \ref{wowgen}, neither do we rule it out.)
The trilinear endpoint Kakeya situation in three dimensions corresponds
to $n=d=3$, $p_1=p_2=p_3 =1/2$ and the matrices $A_j$ each having rank
$2$. Interestingly, in this case $Q_{p}(t)$ is \emph{not} monotonic
decreasing, as the following lemma will demonstrate.
\begin{lemma}\label{notmonlemma}
Let $d \geq 3$. For $j = 1,2, \hdots, d$ suppose $A_j$ is a nonnegative
definite $d \times d$ matrix of rank $(d-1)$ such that $\ker A_1, \ker A_2,
\hdots \ker A_d$ span. Suppose furthermore that $\tfrac{1}{d-1}(A_1+A_2+
\cdots + A_d)\geq A_1,A_2, \hdots ,A_d$. Then $A_d$ is uniquely determined
by $A_1, A_2, \dots A_{d-1}$.
\end{lemma}
\begin{proof}
According to \cite[Proposition 3.6]{BCCT}, under the hypotheses of
the lemma there exists an invertible $d \times d$ matrix $D$ and rank
$(d-1)$ projections $P_j$ such that $A_j = D^* P_j D$ for all $j$ and
$\tfrac{1}{d-1}(P_1+P_2+ \dots + P_d) = I_d$. We claim that for $j \neq k$,
$\ker P_j$ and $\ker P_k$ are orthogonal. To see this let $x_j$ be a
unit vector in $\ker P_j$
and note that $x_1 = \tfrac{1}{d-1}(P_1+P_2+ \cdots + P_d)x_1 =
\tfrac{1}{d-1}(P_2+ \cdots + P_d)x_1$. So
$\langle \tfrac{1}{d-1}(P_2+ \cdots + P_d)x_1, x_1 \rangle = 1$.
On the other hand $\langle P_j x_1, x_1 \rangle \leq 1$ for all $j$ since
$P_j$ is a projection. Thus $\langle P_j x_1, x_1 \rangle = 1$, and $x_1 =
P_j x_1$ for all $j \neq 1$. Similarly $P_j x_i = x_i$ for all $i \neq j$.
Hence for $i \neq j$, $\langle x_i, x_j \rangle = \langle P_j x_i,
x_j \rangle = \langle x_i , P_j x_j \rangle = 0$. Thus the kernels of the
$P_j$'s are mutually orthogonal. By choosing a suitable orthonormal basis
(and possibly changing $D$) we may assume therefore that $P_j = A_j^0$, the
$j^\th$ Loomis--Whitney matrix.

Now let us suppose that $\tilde{D}$ is another $d \times d$ invertible matrix
and that $\tilde{A}_j = \tilde{D}^* A_j^0 \tilde{D}$. Then the statement
$A_i = \tilde{A}_i$ is the same as $A_i^0 = R^* A_i^0 R$ with $R = \tilde{D}
D^{-1}$; i.e. $R$ leaves both $\ker A^0_i$ and its orthogonal complement
$\im A_i^0$ invariant and acts as an isometry on the latter. If now
$A_1 = \tilde{A}_1, \dots, A_{d-1} = \tilde {A}_{d-1}$, this means that
$R$ acts as an isometry on every coordinate hyperplane except possibly
$\{x_d = 0\}$, and the standard basis vectors $e_1, \dots , e_{d-1}$ are
eigenvectors of $R$. This forces the matrix of $R$ with respect to the
standard basis to be diagonal with entries $\pm 1$ and thus
$A_d = \tilde{A}_d$ too.
\end{proof}
\begin{proposition}\label{notmon}
For $j=1,2, \dots , d$ let $W_j$ be a set of nonnegative definite
rank $(d-1)$ $d\times d$ matrices such that for all $A_j\in W_j$,
$\ker A_1, \ker A_2, \hdots ,\ker A_d$ span. Let $\mu_j$ be a positive
Borel measure on $W_j\times\R^d$. Let
$p=(1/(d-1),\hdots,1/(d-1))$.
If $Q_p(1)\leq Q_p(0)$ for all such
positive measures $\mu_j$, then either $d=2$ or each $W_j$ is a singleton,
and $\tfrac{1}{d-1}(A_1+A_2+ \cdots + A_d)\geq A_1, \dots , A_d$.
\end{proposition}
\begin{remark}
In the latter case $d \geq 3$ of Proposition \ref{notmon} the family
$\{A_j\}$ is an affine image of the Loomis--Whitney matrices $\{A_j^0\}$,
i.e. for some invertible $D$, $\tilde{A}_j = \tilde{D}^* A_j^0 \tilde{D}$,
as the proof of Lemma \ref{notmonlemma} shows.
\end{remark}
\begin{proof}
Let $\mu_j^{\#}$ be arbitrary positive finite Borel measures on $\R^d$, and
let $A_j \in W_j$. With
$$f_{j}(t,x):=\int_{{\bf R}^d}
e^{-\pi \langle A_{j}(x-v_jt), (x-v_jt) \rangle}\ d\mu_{j}^{\#}(v_j)$$
and
$$Q^{\#}_p(t):=\int_{{\bf R}^d}f_{1}(t,x)^{1/(d-1)} \cdots
f_{d}(t,x)^{1/(d-1)}dx, $$
setting $\mu_j:=\mu_j^{\#}\otimes\delta_{0}(A_j)$
we now have $Q^\# _p (1) \leq Q^\# _p (0)$. By \cite[Proposition 3.6]{BCCT}
this forces
$$\tfrac{1}{d-1}(A_1+A_2+ \cdots + A_d)\geq A_1, \dots, A_d.$$
When $d \geq 3$, Lemma \ref{notmonlemma} shows that any $d-1$ of the
$A_j$'s determine the remaining one. Thus no $W_j$ may contain more
than one point.
\end{proof}

\section{Lower levels of multilinearity}

As we remarked in the introduction, if $n<d$ then the
conjectured exponents for $n$-linear restriction type problems
depend on the curvature properties of submanifolds in question; with flatter
surfaces expected to enjoy fewer restriction estimates.  Our methods here cannot easily
take advantage of such curvature hypotheses, though, and we will instead address
the issue of establishing $n$-linear restriction estimates in $\R^d$ which assume only
transversality properties rather than curvature properties.  In such a case we can
establish quite sharp estimates.

A similar situation occurs when considering $n$-linear Kakeya estimates in $\R^d$.  The analogue of
``curvature'' would be some sort of direction separation condition (or perhaps a mixed Lebesgue
norm condition) on the tubes in a given family.  Here we will only consider estimates
in which the tubes are counted by cardinality rather than in mixed norms,
and transversality is assumed rather than direction separation.

We begin by discussing the $n$-linear Kakeya situation.
Fix $3 \leq n \leq d$ (the case $n=1$ turns out to be void, and the $n=2$ case standard, see e.g. \cite{TVV}).
Suppose $\tubes_{1},\hdots,\tubes_{n}$ are families of
$\delta$-tubes in $\R^d$. Suppose further that for each $1\leq
j\leq n$, the tubes in $\tubes_{j}$ have long sides pointing in
directions belonging to some sufficiently small \textit{fixed}
neighbourhood of the $j$th standard basis vector $e_{j}$ in ${\bf
S}^{d-1}$.
(Again, the vectors $e_{1},\hdots,e_{n}$ may be replaced by any
fixed set of $n$ linearly independent vectors in $\R^d$ here, as
affine invariance considerations reveal.)
\begin{theorem}\label{nlinrd}
If $\tfrac{n}{n-1}<q\leq\infty$ then there exists a constant $C$,
independent of $\delta$ and the families of tubes
$\tubes_{1},\hdots,\tubes_{n}$, such that
\begin{equation}\label{mkn}
\Bigl\|\prod_{j=1}^{n}\Bigl(\sum_{T_{j}\in\tubes_{j}}\chi_{T_{j}}\Bigr)\Bigr\|_{L^{q/n}({\bf
R}^d)}\leq C\prod_{j=1}^{n}(\delta^{d/q}\: \mbox{\#}\tubes_{j}).
\end{equation}
\end{theorem}
\begin{remark} One can conjecture the same result to hold at the endpoint $q=\tfrac{n}{n-1}$; note for instance
that the estimate is easily verifiable at this endpoint when $n=2$.
For $q < \tfrac{n}{n-1}$ the estimate is false, as can be seen by choosing each family
$\tubes_{j}$ to be (essentially) a partition of a $\delta$-neighbourhood of the
unit cube in $\mbox{span}\{e_{1},\hdots,e_{n}\}\subset\R^{d}$ into parallel tubes oriented in the direction $e_j$.
It is likely that the estimates can be improved if some direction separation condition is imposed on each of the $\tubes_j$, but
we do not pursue this matter here.
\end{remark}

In order to prove Theorem \ref{nlinrd} we first apply a rescaling
(in the spirit of \eqref{1tubes}) to reduce inequality \eqref{mkn}
to an equivalent statement for families of tubes of width $1$, and
arbitrary length. Then we simply dominate the characteristic
functions of these dilated tubes by appropriate gaussians, and
appeal to Corollary \ref{corwow}, setting $p_j=p>\tfrac{1}{n-1}$
and $M_j=A_{j}^{0}$, where $A_{j}^{0}$ is the $j^\th$
Loomis--Whitney matrix in $\R^d$, and $1\leq j\leq n$.

The corresponding $n$-linear restriction inequalities (with the
familiar $\eps$-loss in the localisation parameter $R$) may
now be obtained by the bootstrapping argument from Section 2, thus obtaining an estimate of the form
\begin{equation}\label{rest-loss}
\Bigl\|\prod_{j=1}^{n}\ex_{j}g_{j}\Bigr\|_{L^{q/n}(B(0,R))}\leq C
R^{\eps}\prod_{j=1}^{n} \|g_{j}\|_{L^p(U_j)}
\end{equation}
when $q \geq \frac{2n}{n-1}$ and $p' \leq \frac{n-1}{n} q$.
What
is perhaps particularly curious is the fact that for $n<d$, the
standard Rademacher-function argument \textit{does not} allow the
optimal $n$-linear Kakeya inequalities to be obtained from the
corresponding optimal $n$-linear restriction inequalities.
To avoid
repetition of the arguments in Section 2, we omit the details.

\begin{remark}
The epsilon loss in \eqref{rest-loss} should be removable (thus allowing $R$ to be sent to infinity);
certainly this is possible in the $n=2$ case, see for instance \cite{TaoRestn} for this standard and useful
estimate.  The conditions
$q \geq \frac{2n}{n-1}$ and $p' \leq \frac{n-1}{n} q$ can be verified to be sharp (e.g. by considering the
Loomis--Whitney case when the maps $\Sigma_j$ are linear) but can be improved when the $\Sigma_j$ have additional
curvature properties (again, see for instance \cite{TaoRestn} for a survey).
\end{remark}

\begin{remark}
One may also obtain non-trivial multilinear estimates for
$k$-plane transforms from Corollary \ref{corwow} by choosing the
matrices $M_j$ to be appropriate projections onto
$(d-k)$-dimensional subspaces of $\R^d$. For example, if
$M_{j}=I-A_{j}^{0}$ and $n=d$, Proposition \ref{wowgen} implies
certain multilinear analogues of the Radon transform estimates of Oberlin
and Stein \cite{OS}.\footnote{It turns out that when all the
$A_j$'s have rank $1$, a non-perturbative linear analysis is
rather straightforward. We shall return to such matters at a later
date.} We leave the details of these implications to the
interested reader.
\end{remark}

\section{Variable-coefficient extensions}\label{variable-sec}
More general (diffeomorphism-invariant) families of oscillatory
integral operators, of which the extension operators are examples,
were first considered by H\"ormander in \cite{Hor}. H\"ormander
conjectured that under certain natural non-degeneracy conditions
on the associated phase function (see \cite{St2}), such operators
would satisfy $L^{p} \to L^{q}$ estimates in agreement with the
classical restriction conjecture. It is now well-known that this
conjecture is in general false -- see Bourgain \cite{Bofalse}. In
this section we consider the validity of such
generalisations of the
\textit{multilinear} restriction problem discussed in the
introduction, and obtain almost optimal results in this setting.

Let $\Phi :\R^{d-1}\times\R^{d}\rightarrow\R$ be a smooth phase
function, $\lambda>0$ and $\psi
:\R^{d-1}\times\R^{d}\rightarrow\R$ be a compactly supported
smooth cut-off function. We define the operator $S_{\lambda}$ by
$$
S_{\lambda}g(\xi):=\int_{{\bf
R}^{d-1}}e^{i\lambda\Phi(x,\xi)}\psi(x,\xi)g(x)dx,
$$
and the vector field $X(\Phi)$ by
$$
X(\Phi):=\bigwedge_{k=1}^{d-1}\frac{\partial}{\partial
x_{k}}\nabla_{\xi}\Phi.
$$

\begin{remark} When the phase $\Phi$ takes the form $\Phi(x,\xi) = x \cdot \Sigma(\xi)$
then the operator $S_\lambda$ is essentially (up to rescaling and cutoffs) an extension
operator $\ex$, and $X(\Phi)$ is essentially the vector field $Y$.
\end{remark}

Now we suppose that $S_{\lambda}^{(1)},\hdots, S_{\lambda}^{(d)}$
are such operators associated to phase functions
$\Phi_{1},\hdots,\Phi_{d}$, and cut-off functions
$\psi_{1},\hdots,\psi_{d}$.

Our generalisation of the multilinear transversality condition,
which we will impose from this point on, will be that for some
constant $\nu>0$,
\begin{equation}\label{difftrans}
\det\left(X(\Phi_{1})(x^{(1)},\xi),\hdots,X(\Phi_{d})(x^{(d)},\xi)\right)>\nu
\end{equation}
for all
$(x^{(1)},\xi)\in\supp(\psi_{1}),\hdots,(x^{(d)},\xi)\in\supp(\psi_{d})$.
In addition to this, for each multi-index $\beta\in N^{d-1}$ let
us suppose that for some constant $A_{\beta}\geq 0$,
\begin{equation}\label{sm}
\|\partial_{x}^{\beta}\Phi_{j}(x,\cdot)\|_{C^{2}_{\xi}({\bf R}^d)}\leq
A_{\beta}\;\;\;\mbox{for all } 1\leq j\leq d,\; x \in \R^{d-1}.
\end{equation}
\begin{theorem}\label{diffrest}
If \eqref{difftrans} and \eqref{sm} hold, then for each
$\eps>0$, $q\geq\tfrac{2d}{d-1}$ and $p'\leq\tfrac{d-1}{d}q$,
there is a constant $C>0$, depending only on $\eps$, $p$, $q$,
$d$, $\nu$ and finitely many of the $A_{\beta}$'s, for which
$$
\Bigl\|\prod_{j=1}^{d}S_{\lambda}^{(j)}g_{j}\Bigl\|_{L^{q/d}({\bf
R}^{d})} \leq
C_{\eps}\lambda^{\eps}\prod_{j=1}^{d}\lambda^{-d/q}\|g_{j}\|_{L^{p}({\bf
R}^{d-1})}
$$
for all $g_{1},\hdots,g_{d}\in L^{p}({\bf R}^{d-1})$ and $\lambda>
0$.
\end{theorem}
\begin{remark}\label{exglw}
It is possible that the above inequality continues to hold with
$\eps=0$, although we have been unable to prove this; for instance, the number of $A_\beta$ which
we need in our argument goes to infinity as $\eps \to 0$, though this may well be unnecessary.
It has been shown recently (in \cite{BCW}, by very different
techniques), that if $\Phi_{j}(x,\xi)=x\cdot\Gamma_{j}(\xi)$,
where $\Gamma_{j}:\R^{d}\rightarrow\R^{d-1}$ are smooth
submersions, then one may indeed set $\eps=0$ in the
conclusion of Theorem \ref{diffrest}. In this particular example,
$S_{\lambda}^{(j)}g_{j}(\xi)$ is essentially
$\widehat{g}_{j}\circ\Gamma_{j}(\lambda\xi)$, and so at the sharp
endpoint ($p=2$, $q=\tfrac{2d}{d-1}$), by Plancherel's theorem
one may see this example as a non-linear generalisation of the
Loomis--Whitney inequality \eqref{clw}.
\end{remark}
\begin{remark}
The implicit conditions $X(\Phi_{1}),\hdots, X(\Phi_{d})\not =0$
are of course equivalent to the statement that the matrices
$\tfrac{\partial^{2}\Phi_{1}}{\partial
x\partial\xi},\hdots,\tfrac{\partial^{2}\Phi_{d}}{\partial
x\partial\xi}$ all have full rank $d-1$. We point out that in
general, given \eqref{difftrans}, one cannot expect further
non-degeneracy assumptions on the phase functions $\Phi_{j}$ to
lead to improvements in the claimed range of exponents here.
Again, this is very much in contrast with what happens at lower
levels of multilinearity.
\end{remark}

We now come to the corresponding variable-coefficient multilinear
Kakeya-type problem. For a discussion of the original linear
setting see Wisewell \cite{wisewell}.

For each $1\leq j\leq d$ let $\tubes_{j}$ denote a collection of
subsets of $\R^{d}$ of the form
$$\{\xi\in\R^{d}:
|\nabla_{x}\Phi_{j}(a,\xi)-\omega|\leq\delta,\;
(a,\xi)\in\supp(\psi_j)\},
$$
where $a,\omega\in\R^{d-1}$. It is important for us to observe
that the conditions \eqref{difftrans} and \eqref{sm} (with
$|\beta| \leq 1$) imply that these sets contain, and are contained in,
$O(\delta)$-neighbourhoods of smooth curves in $\R^d$. For this
reason it is convenient to extend the use of our tube notation and
terminology from the previous sections. The implicit constants in
the $O$-notation here depend on $d$, $\nu$ and the $A_{\beta}$'s
with $|\beta| \leq 1$ (these quantities appear in \eqref{sm}).

\begin{theorem}\label{kakdiff} If \eqref{difftrans} and \eqref{sm} hold, then
for each $\eps>0$ and $q\geq \tfrac{d}{d-1}$ there exists a
constant $C>0$, depending only on $\eps$, $q$, $d$, $\nu$ and
finitely \footnote{In fact one only really sees the $A_{\beta}$'s
with $|\beta|=1$ here -- this is easily seen from our proof.} many
of the $A_{\beta}$'s, such that
$$
\Bigl\|\prod_{j=1}^{d}\sum_{T_{j}\in\tubes_{j}}\chi_{T_{j}}\Bigl\|_{L^{q/d}({\bf
R}^d)}\leq
C_{\eps}\delta^{-\eps}\prod_{j=1}^{d}\left(\delta^{d/q}\mbox{\#}\tubes_{j}\right)
$$
for all collections $\tubes_{1},\hdots,\tubes_{d}$ and $\delta>0$.
\end{theorem}
\begin{remark}
Again, it seems possible that the above inequality continues to hold
with $\eps=0$. Notice that in the special case described in
Remark \ref{exglw}, the corresponding $T_j\in\tubes_j$ are simply
unions of fibres of the submersions $\Gamma_j$.
\end{remark}
The proofs of both Theorems \ref{diffrest} and \ref{kakdiff}
follow bootstrapping arguments closely related to that of Bourgain
used in Section 2.\footnote{Variants of such bootstrapping
arguments have been considered previously by both Wolff and
the third author, although not in a multilinear setting.} For these we
need some further notation:

For $\alpha>0$, $q\geq\tfrac{2d}{d-1}$ and
$p'\leq\tfrac{d-1}{d}q$, let
$$\re^{*}_{c}(p\times\cdot\cdot\cdot\times p\rightarrow
q;\alpha)$$ denote the multilinear oscillatory integral estimate
$$
\Bigl\|\prod_{j=1}^{d}S_{\lambda}^{(j)}g_{j}\Bigl\|_{L^{q/d}({\bf
R}^{d})} \leq
C\lambda^{\alpha}\prod_{j=1}^{d}\lambda^{-d/q}\|g_{j}\|_{L^{2}({\bf
R}^{d-1})},
$$
for some constant $C$ (depending only on $\alpha$, $p$, $q$, $d$,
$\nu$ and finitely many of the $A_{\beta}$'s), all
$g_{1},\hdots,g_{d}\in L^{2}({\bf R}^{d-1})$ and $\lambda> 0$.

Similarly, for each $\alpha>0$ and $q\geq\tfrac{d}{d-1}$ let
$$\ka^{*}_{c}(1\times\cdot\cdot\cdot\times 1\rightarrow q;\alpha)
$$
denote the multilinear ``curvy" Kakeya estimate
\begin{equation}\label{subficialdiff}
\Bigl\|\prod_{j=1}^{d}\sum_{T_{j}\in\tubes_{j}}\chi_{T_{j}}\Bigl\|_{L^{q/d}
({\bf R}^d)}\leq C
\delta^{-\alpha}\prod_{j=1}^{d}\left(\delta^{d/q}\mbox{\#}\tubes_{j}\right),
\end{equation}
for a similarly uniform constant $C$, all $\delta>0$ and all
families $\tubes_{1},\hdots,\tubes_{d}$. (Of course the hypotheses
\eqref{difftrans} and \eqref{sm} are assumed implicitly here.)
\begin{remark}
Inequality \eqref{subficialdiff} is easily seen to be equivalent
to the superficially stronger
\begin{equation}\label{superficialdiff}
\Bigl\|\prod_{j=1}^{d}\sum_{T_{j}\in\tubes_{j}}\lambda_{T_{j}}\chi_{T_{j}}\Bigl\|_{L^{q/d}
({\bf R}^d)}\leq C
\delta^{-\alpha}\prod_{j=1}^{d}\Bigl(\delta^{d/q}\sum_{T_{j}\in\tubes_{j}}\lambda_{T_{j}}\Bigr),
\end{equation}
uniformly in the non-negative constants $\lambda_{T_{j}}$
($T_{j}\in\tubes_{j}$, $1\leq j\leq d$).
\end{remark}
Given that the estimate $\ka^*(1\times\cdot\cdot\cdot\times
1\rightarrow q;0)$ holds for all $q>\tfrac{d}{d-1}$ (Theorem
\ref{kakthm}), we may reduce the proof of Theorem \ref{kakdiff} to
a repeated application of the following. Note that by H\"older's
inequality it suffices to treat $q>\tfrac{d}{d-1}$.
\begin{proposition}\label{prokak}
For each $\alpha,\eps>0$ and $q>\tfrac{d}{d-1}$,
$$
\ka_{c}^{*}(1\times\cdot\cdot\cdot\times 1\rightarrow q;\alpha)+
\ka^{*}(1\times\cdots\times 1\rightarrow q;\eps) \implies
\ka_{c}^{*}(1\times\cdot\cdot\cdot\times 1\rightarrow
q;\tfrac{\alpha}{2}+\tfrac{\eps}{2}).
$$
\end{proposition}
Given Theorem \ref{kakdiff} (and H\"older's inequality), we may
similarly reduce the proof of Theorem \ref{diffrest} to a
bootstrapping argument.\footnote{In principle one ought to be able
to prove Theorem \ref{diffrest} directly using a suitable variant
of Theorem \ref{restthm} along with Theorem \ref{kakdiff}. We do
not pursue this matter here.}
\begin{proposition}\label{prorest}
For each $\alpha,\eps,\eps_{0}>0$ and $q>\tfrac{2d}{d-1}$,
$$
\re_{c}^{*}(2\times\cdot\cdot\cdot\times 2\rightarrow q;\alpha)+
\ka_{c}^{*}(1\times\cdot\cdot\cdot\times 1\rightarrow
\tfrac{q}{2};\eps) \implies
\re_{c}^{*}(2\times\cdot\cdot\cdot\times 2\rightarrow
q;\tfrac{\alpha}{2}+\tfrac{\eps}{4}+\eps_{0}).
$$
\end{proposition}
\subsection*{Proof of Proposition \ref{prokak}} Let $\{B\}$ be a
tiling of $\R^{d}$ by cubes of side $\sqrt{\delta}$, and write
$$
\Bigl\|\prod_{j=1}^{d}\sum_{T_{j}\in\tubes_{j}}\chi_{T_{j}}
\Bigl\|_{L^{q/d}({\bf R}^d)}^{q/d}
=\sum_{B}\Bigl\|\prod_{j=1}^{d}\sum_{T_{j}\in\tubes_{j}}\chi_{T_{j}\cap
B}\Bigl\|_{L^{q/d}({\bf R}^d)}^{q/d}.
$$
Now, by the smoothness of the $\Phi_j$'s, each $T_{j}\cap B$ is
contained in a rectangular tube of dimensions
$O(\delta)\times\cdots\times O(\delta)\times O(\sqrt{\delta})$,
and furthermore (upon rescaling) the $d$ families of these
``rectangular" tubes $\{T_{j}\cap B\}_{T_{j}\in\tubes_{j}}$
($1\leq j\leq d$) have the transversality property required by the
hypothesis $\ka^{*}(1\times\cdot\cdot\cdot\times 1\rightarrow
q;\eps)$. Hence
$$
\Bigl\|\prod_{j=1}^{d}\sum_{T_{j}\in\tubes_{j}}\chi_{T_{j}\cap
B}\Bigl\|_{L^{q/d}({\bf R}^d)} \lesssim
\delta^{-\eps/2}\prod_{j=1}^{d}\delta^{d/q}\mbox{\#}\{
T_{j}\in\tubes_{j}:T_{j}\cap B\not=\emptyset\}
$$
uniformly in $B$. We note that the implicit constants in the $O$
notation above depend only on $d$, and the constants $\nu$ and
$A_{\beta}$ with $|\beta|=1$.\footnote{Here we are using the
$C^{2}_{\xi}$-control of the phase functions $\Phi_{j}$ to
guarantee that the tubes are ``locally straight" in the claimed
way.} We next observe the elementary fact that
$$
\mbox{\#}\{ T_{j}\in\tubes_{j}:T_{j}\cap B\not=\emptyset\}
\lesssim
\sum_{T_{j}\in\tubes_{j}}\chi_{T_{j}+B(0,c\sqrt{\delta})}(\xi_{B}),
$$
uniformly in $\xi_{B}\in B$, where $c$ is a sufficiently large
constant, again depending only on $d$, $\nu$ and $A_{\beta}$ with
$|\beta|= 1$. Now each $T_{j}$ is given by
$$
T_{j}=\{\xi\in\R^{d}:
|\nabla_{x}\Phi_{j}(a_{j},\xi)-\omega_{j}|\leq\delta,\;
(a_{j},\xi)\in\supp(\psi_{j})\}
$$
for some $a_{j},\omega_{j}\in\R^{d-1}$, and so by our geometric
interpretation of such tubes we see that
$
T_{j}+B(0,c\sqrt{\delta})
\subset\widetilde{T}_{j}$, where
$$\widetilde{T}_{j}:=\left\{\xi\in\R^d:|\nabla_{x}\Phi_{j}(a_j,\xi)-\omega_j|\lesssim\sqrt{\delta},\;
(a_{j},\xi)\in\supp(\psi_{j}) + B(0, O(\sqrt{\delta}))\right\},
$$
yielding
\begin{eqnarray*}
\begin{aligned}
\Bigl\|\prod_{j=1}^{d}\sum_{T_{j}\in\tubes_{j}}\chi_{T_{j}}\Bigl\|_{L^{q/d}({\bf
R}^d)}
\lesssim\delta^{-\eps/2}\biggl(\sum_{B}\Bigl(\prod_{j=1}^{d}\delta^{d/q}\sum_{\widetilde{T}_{j}}
\chi_{\widetilde{T}_{j}}(\xi_{B})\Bigr)^{q/d}\biggr)^{d/q}
\end{aligned}
\end{eqnarray*}
uniformly in the choice of $\xi_{B}$. Hence upon averaging we
obtain
\begin{eqnarray*}
\begin{aligned}
\Bigl\|\prod_{j=1}^{d}\sum_{T_{j}\in\tubes_{j}}\chi_{T_{j}}\Bigl\|_{L^{q/d}({\bf
R}^d)}\lesssim \delta^{-\eps/2+d^{2}/(2q)}
\Bigl\|\prod_{j=1}^{d}\sum_{\widetilde{T}_{j}:T_{j}\in\tubes_{j}}
\chi_{\widetilde{T}_{j}}\Bigl\|_{L^{q/d}({\bf R}^d)},
\end{aligned}
\end{eqnarray*}
which by the remaining hypothesis\footnote{There is an extremely minor issue here,
which is that the tubes $\tilde T_j$ are defined with $(a_j,\xi)$ ranging in
$\supp(\psi_{j}) + B(0, O(\sqrt{\delta}))$ rather than $\supp(\psi_{j})$.  But this negligible enlargement of the support can be dealt with
by modifying $\psi_j$ very slightly and checking that the various bounds on the geometry do not change very much.  We omit the details.}
 $\ka^{*}_{c}(1\times\cdots\times
1\times q;\alpha)$ is
$$
\lesssim\delta^{-\eps/2-\alpha/2}\prod_{j=1}^{d}\delta^{d/q}\mbox{\#}\tubes_{j},
$$
completing the proof of the proposition.

\begin{remark}
Our approach to Proposition \ref{prorest} is somewhat different
technically from that of Proposition \ref{MRiffMK}. This is
largely due to our desire to avoid formulating a potentially
cumbersome generalisation of Lemma \ref{uncertainty}. The drawback
of the resulting (slightly cruder) argument, which is largely
aesthetic, is the additional epsilon-loss present in the
statement, and the role played by the high order derivatives in
\eqref{sm}.
\end{remark}

\subsection*{Proof of Proposition \ref{prorest}} The argument we give
here has essentially one additional ingredient to that of
Proposition \ref{prokak} -- a ``wavepacket" decomposition, and
thus one is forced to deal with the additional technicalities
associated with the uncertainty principle; i.e. the fact that
wavepackets are not genuinely supported on tubes, but rather decay
rapidly away from them. Let $\{Q\}$ be a tiling of $\R^{d-1}$ by
cubes $Q$ of side $\lambda^{-1/2}$, and for such a $Q$, let $x_{Q}$ be
its centre. We now decompose each $g_{j}$
into local Fourier series at an appropriate scale. Let
$\{\chi_{Q}\}$ be a (smooth) partition of unity adapted to the
tiling $\{Q\}$. For uniformity purposes let us suppose that
$\chi_{Q}(x)=\chi(\lambda^{1/2}(x-x_Q))$ for some smooth compactly
supported function $\chi$. Now, for each $g_{j}$ we may write
$$
g_{j}=\sum_{Q} \sum_{\ell\in\lambda^{-1/2}{\bf
Z}^{d-1}}a_{Q,\ell}^{(j)}e_{Q,\ell},
$$
where $e_{Q,\ell}$ is the modulated cap
$$e_{Q,\ell}(x):=\chi_{Q}(x)e^{i\ell\cdot x},
$$
and the $a_{Q,\ell}^{(j)}$'s are complex numbers. By linearity of
$S_{\lambda}^{(j)}$,
$$
S_{\lambda}^{(j)}g_{j}=\sum_{Q,\ell}a_{Q,\ell}^{(j)}S_{\lambda}^{(j)}e_{Q,\ell}.
$$

Now we localize the $S_{\lambda}^{(j)}e_{Q,\ell}$ to tubes.
Let $\eta := q\eps_{0}/d^{2}$. For each $\ell\in
\lambda^{-1/2}\Z^{d-1}$ and $1\leq j\leq d$ let
$R^{(j)}_{Q,\ell}$ be the curved tube
$$
R^{(j)}_{Q,\ell}:=\{\xi\in\R^{d}:
|\nabla_{x}\Phi_{j}(x_{Q},\xi)-\ell|\leq\lambda^{-1/2+\eta},\;
(x_Q,\xi)\in\supp(\psi_j)\}.
$$
By a standard repeated integration by parts argument we have that
for each $M\in\N$,
\begin{equation}\label{princ}
|S_{\lambda}^{(j)}e_{Q,\ell}(\xi)|\lesssim\lambda^{-M\eta-(d-1)/2},
\end{equation}
for all $\xi\in \R^{d}\backslash R_{Q,\ell}^{(j)}$. Naturally the
implicit constants here depend on $d$ and the smoothness bounds
$A_{\beta}$ for $|\beta|\leq M$.

We now tile $\R^{d}$ by cubes $B$ of side $\lambda^{-1/2}$. The
idea is to use the oscillatory integral estimate in our hypothesis
on each $B$, and then use the curvy Kakeya estimate to reassemble
them. Let $\mathcal{P}$ denote the set of non-empty subsets $P$ of
the set of integers $\{1,\hdots,d\}$, and for each
$P\in\mathcal{P}$ let $P^c$ denote the complement of $P$ in
$\{1,\hdots, d\}$. Now by the triangle inequality
$$
\Bigl\|\prod_{j=1}^{d}S_{\lambda}^{(j)}g_{j}\Bigr\|_{L^{q/d}({\bf
R}^{d})}=
\Biggl(\sum_{B}\Bigl\|\prod_{j=1}^{d}S_{\lambda}^{(j)}g_{j}
\Bigr\|_{L^{q/d}(B)}^{q/d}\Biggr)^{d/q}\leq I
+\sum_{P\in\mathcal{P}}I_{P},
$$
where
$$
I:= \Biggl(\sum_{B}\Bigl\|\prod_{j=1}^{d}\Bigl(\sum_{Q,\ell:
R_{Q,\ell}^{(j)}\cap
B\not=\emptyset}a_{Q,\ell}^{(j)}S_{\lambda}^{(j)}e_{Q,\ell}\Bigr)\Bigr\|
_{L^{q/d}(B)}^{q/d}\Biggr)^{d/q}
$$
and $I_{P}$ is given by $I_{P}^{q/d}:=$
$$
\sum_{B}\Bigl\|\Bigl(\prod_{j\in P}
\sum_{Q,\ell:R_{Q,\ell}^{(j)}\cap B=\emptyset}
a_{Q,\ell}^{(j)}S_{\lambda}^{(j)}e_{Q,\ell}\Bigr)\Bigl(\prod_{k\in
P^c} \sum_{Q,\ell:R_{Q,\ell}^{(k)}\cap B\not=\emptyset}
a_{Q,\ell}^{(k)}S_{\lambda}^{(k)}e_{Q,\ell}\Bigr)
\Bigr\|_{L^{q/d}(B)}^{q/d}.
$$
We first estimate the principal term $I$. By rescaling the
hypothesis $\re^{*}_{c}(2\times\cdots\times 2\rightarrow
q;\alpha)$, and observing the scale-invariance\footnote{
There is a minor technical issue here.
The terms in condition \eqref{sm} containing
zero $\xi$-derivatives actually fail to be
invariant in the appropriate way. However, this may be easily
rectified by subtracting off harmless affine factors from the phases
$\Phi_j$, and absorbing
them into the functions $g_j$.}
of conditions
\eqref{difftrans} and \eqref{sm}, we have that
$$I\lesssim \lambda^{\alpha/2+d^{2}/q}
\Biggl(\sum_{B}\prod_{j=1}^{d}\Bigl\|\sum_{Q,\ell:
R_{Q,\ell}^{(j)}\cap
B\not=\emptyset}a_{Q,\ell}^{(j)}e_{Q,\ell}\Bigr\|
_{2}^{q/d}\Biggr)^{d/q},
$$
which by the almost orthogonality of the $e_{Q,\ell}$'s is further
bounded by
\begin{eqnarray*}
\begin{aligned}
\lambda^{\alpha/2-d^{2}/q}
\Biggl(&\sum_{B}\prod_{j=1}^{d}\Bigl(\sum_{Q}\lambda^{-(d-1)/2}
\sum_{\ell: R_{Q,\ell}^{(j)}\cap
B\not=\emptyset}|a_{Q,\ell}^{(j)}|^{2} \Bigr)
^{q/(2d)}\Biggr)^{d/q}\\
&\lesssim\lambda^{\alpha/2-d^{2}/q-d(d-1)/4}
\Biggl(\sum_{B}\prod_{j=1}^{d}\Bigl(\sum_{Q,\ell}
|a_{Q,\ell}^{(j)}|^{2}\chi_{R_{Q,\ell}^{(j)}}(\xi_{B}) \Bigr)
^{q/(2d)}\Biggr)^{d/q},
\end{aligned}
\end{eqnarray*}
uniformly in $\xi_{B}\in B$. Strictly speaking the tubes in this
last line should be replaced by slightly dilated versions of
themselves, however we shall gloss over this detail. On averaging
and applying the remaining hypothesis
$\ka^*_{c}(1\times\cdots\times 1\rightarrow\tfrac{q}{2};\eps)$
(in the equivalent form \eqref{superficialdiff}) to the above
expression we obtain
\begin{eqnarray*}
\begin{aligned}
I&\lesssim\lambda^{\alpha/2-d^{2}/q-d(d-1)/4}
\Biggl(\lambda^{-d/2}\int_{{\bf
R}^d}\Bigl(\prod_{j=1}^{d}\sum_{Q,\ell}
|a_{Q,\ell}^{(j)}|^{2}\chi_{R_{Q,\ell}^{(j)}}(x) \Bigr)
^{q/(2d)}dx\Biggr)^{d/q}\\
&\lesssim\lambda^{\alpha/2-d^{2}/(2q)-d(d-1)/4+\eps(1-2\eta)/4-(1-2\eta)d^{2}/(2q)}
\Bigl(\prod_{j=1}^{d}\sum_{Q,\ell}|a_{Q,\ell}^{(j)}|^{2}\Bigr)^{1/2}\\
&\lesssim
\lambda^{\alpha/2+\eps/4+\eps_{0}}\prod_{j=1}^{d}\lambda^{-d/q}\|g_{j}\|_{2}.
\end{aligned}
\end{eqnarray*}
Here in this last line we have used Plancherel's theorem and the
fact that $\eta=q\eps_{0}/d^{2}$.

It is enough now to show that for each $P\in\mathcal{P}$ and $N>0$
we have the error estimates
$$I_{P}\lesssim\lambda^{-N\eta}\prod_{j=1}^{d}\|g_{j}\|_{2}.$$
However, this is an elementary consequence of H\"older's
inequality and the decay estimate \eqref{princ}.  This proves Proposition \ref{prorest}.

\section{An application to the joints problem}\label{joints-sec}

We now give an application of the multilinear Kakeya estimate
(Theorem \ref{kakthm}) to a discrete geometry problem, namely the
``joints'' problem studied in \cite{chaz}, \cite{sharir},
\cite{feldman}.

Let us recall the setup for this problem.  Consider a collection
$L$ of $n$ lines in $\R^3$.  Define a \emph{joint} to be a point
in $\R^3$ which is contained in at least one triple $(l,l',l'')$
of concurrent lines in $L$ which are not coplanar.  (Note that a
single joint may arise from multiple triples, but in such cases we
only count those joints once.) The joints problem is to determine,
for each fixed $n$, the maximum number of joints one can attain
for a configuration $L$ of $n$ lines. This problem was observed to
be formally related to the Kakeya problem in \cite{wolff:survey};
in this paper we establish for the first time a rigorous
connection between the two problems.

An easy lattice construction (where the lines are parallel to the
co-ordinate axes and have two of the co-ordinates fixed to be
integers between $1$ and $\sqrt{n}$) shows that one can have at
least $\sim n^{3/2}$ joints.  In the other direction, one
trivially observes that each line in $L$ can contain at most $n$
joints, and hence we have an upper bound of $n^2$ for the total
number of joints.  There has been some progress in improving the
upper bound; the most recent result in \cite{chaz} shows that the
number of joints is at most $O( n^{112/69} \log^{6/23} n ) \leq O(
n^{1.6232} )$.  It is tentatively conjectured that the lower bound
of $n^{3/2}$ is essentially sharp up to logarithms.

It turns out that the multilinear Kakeya estimate in Theorem
\ref{kakthm} can support this conjecture, provided that the joints
are sufficiently transverse.  For any $0 < \theta \leq 1$, let us
say that three concurrent lines $l, l', l''$ are
\emph{$\theta$-transverse} if the parallelopiped generated by the
unit vectors parallel to $l, l', l''$ (henceforth referred to as the
\emph{directions} of $l$, $l'$, $l''$) has volume at least
$\theta$.  Let us define a \emph{$\theta$-transverse joint} to be
a point in $\R^3$ which is contained in at least one triple
$(l,l',l'')$ of $\theta$-transverse concurrent lines in $L$.  Note
that every joint is $\theta$-transverse for some $\theta$.

\begin{theorem}\label{joints-thm}  For any $0 < \theta \leq 1$, the number of
$\theta$-transverse joints is $\;\;$ $O_\eps( n^{3/2+\eps}
\theta^{-1/2-\eps} )$
for any $\eps > 0$, where the subscripting of $O$ by $\eps$ means
that the implied constant can depend on $\eps$.
\end{theorem}

This theorem suggests that the hard case of the joints problem
arises when considering nearly-coplanar joints, with different
joints being approximately coplanar in different orientations.
This resembles the experience in \cite{KLT}, when the ``plany''
case of the Kakeya problem was by far the most difficult to
handle.

\begin{proof} We first establish this conjecture in the case $\theta \sim 1$.
We cover the unit sphere ${\bf S}^2$ by $O(1)$ finitely
overlapping caps $C_1,\ldots,C_k$ of width $\theta/1000$.  Observe
that if $l,l',l''$ are a $\theta$-transverse collection of lines,
then the directions of $l,l',l''$ will lie in three distinct
caps $C_i, C_{i'}, C_{i''}$, which are transverse in the sense of
\eqref{spanning}. Since the number of such triples of caps is
$O(1)$, it thus suffices to show that
\begin{equation}\label{triple-joint}
\# \{ p \in \R^3: p \in l, l', l'' \hbox{ for some } l \in L_i, l'
\in L_{i'}, l'' \in L_{i''} \} = O(n^{3/2})
\end{equation}
for each such transverse triple $(C_i, C_{i'}, C_{i''})$, where
$L_i$ is the collection of lines in $L$ with directions in
$C_i$.

By rescaling we may assume that all the joints are contained in
the ball of radius $1/1000$ centred at the origin. Let $\delta >
0$ be a small parameter (eventually it will go to zero), and for
each line $l \in L$ let $T_l$ denote the $\delta \times \ldots
\times \delta \times 1$ tube with axis $l$ and centre equal to the
closest point of $l$ to the origin.  Let $\tubes_i$ denote the
collection of all the tubes $T_l$ associated to lines $l$ in
$L_i$, and similarly define $\tubes_{i'}$, $\tubes_{i''}$.  From
elementary geometry we see that if $p$ is an element of the set in
\eqref{triple-joint}, then we have
$$ \sum_{T_i \in \tubes_i} \chi_{T_i}(x) \geq 1$$
whenever $|x-p| < c \delta$, where $c > 0$ is a small absolute
constant depending on the transversality constant of $(C_i,
C_{i'}, C_{i''})$.  Similarly for $\tubes_{i'}$ and
$\tubes_{i''}$.  Since the number of joints is finite, we see that
for $\delta$ sufficiently small, the balls $\{ x \in \R^3: |x-p| <
c\delta \}$ will be disjoint.  We conclude that
$$ \Bigl\| \bigl(\sum_{T_i \in \tubes_i} \chi_{T_{i}}\bigr)
\bigl(\sum_{T_{i'} \in \tubes_{i'}} \chi_{T_{i'}}\bigr)
\bigl(\sum_{T_{i''} \in \tubes_{i''}} \chi_{T_{i''}}\bigr)
\Bigr\|_{L^{q/3}({\bf R}^3)} \geq c_q N^{3/q} \delta^{9/q}$$ for
any $\frac{3}{2} < q \leq \infty$, where $N$ denotes the left-hand
side of \eqref{triple-joint} and $c_q > 0$ is a constant depending
only on $c$ and $q$.  Applying Theorem \ref{kakthm} we obtain
$$ N^{3/q} \delta^{9/q} \leq C_q (\delta^{3/q} \# L_i) (\delta^{3/q} \# L_{i'}) (\delta^{3/q} \# L_{i''})$$
or in other words
\begin{equation}\label{n-bound}
 N \leq C_q^{q/3} (\# L_i \# L_{i'} \# L_{i''})^{q/3}.
\end{equation}
Since $\# L_i, \# L_{i'}, \# L_{i''} \leq n$ and $q$ can be
arbitrarily close to $3/2$, the claim follows.

Now we handle the case when $\theta$ is much smaller than $1$,
using some (slightly inefficient) trilinear variants of the
bilinear rescaling arguments employed in \cite{TVV}.

Suppose that $(l,l',l'')$ are $\theta$-transverse.  Each pair of
lines in $l,l',l''$ determines an angle; without loss of
generality we may take $l, l'$ to subtend the largest angle.
Calling this angle $\alpha$, we see from elementary geometry that
$\theta^{1/2} \lesssim \alpha \lesssim 1$, and that $l''$ makes an
angle of at least $\gtrsim \theta/\alpha$ and at most $\alpha$
with respect to the plane spanned by $l$ and $l'$.  To exploit
this, let us say that $(l,l',l'')$ are
\emph{$(\alpha,\beta)$-transverse} for some $0 < \beta \lesssim
\alpha \lesssim 1$ if $l,l'$ make an angle of $\sim \alpha$ and
$l''$ makes an angle of $\sim \beta$ with respect to the plane
spanned by $l$ and $l'$.  Define a
\emph{$(\alpha,\beta)$-transverse} joint similarly. A simple
dyadic decomposition argument (giving up some harmless factors of
$\log \frac{1}{\theta}$) then show that it suffices to show that
the number of $(\alpha,\beta)$-transverse joints is $O_\eps(
n^{3/2+\eps} (\alpha \beta)^{-1/2+\eps} )$ for every $\eps > 0$.
In fact we will prove the sharper bound of $O_\eps( n^{3/2+\eps}
(\beta/\alpha)^{-1/2+\eps} )$.

Let us first handle the case when $\alpha \sim 1$, so that $l$ and
$l'$ make an angle of $\sim 1$.  By symmetry we may also assume
that $l''$ makes a smaller angle with $l'$ than it does with $l$,
so $l$ and $l''$ also make an angle of $\sim 1$.  By a
decomposition of the sphere into $O(1)$ pieces, we can then assume
that there exist transverse subsets $S, S'$ of the sphere such
that the direction of $l$ lies in $S$, and the directions of
$l'$ and $l''$ lie in $S'$.  (Note that $S'$ may be somewhat
larger than $S$.)

Let $\omega_1, \ldots, \omega_K$ be a maximal $\beta$-separated
set of directions on the sphere, thus $K = O(1/\beta^2)$.  For
each direction $\omega_k$, let $L_k$ denote the family of lines
with direction in $S$ which make an angle of $\pi/2 - O(\beta)$ with
$\omega_k$, thus they are nearly orthogonal to $\omega_k$.  Define
$L'_k$ similarly but with $S$ replaced by $S'$. From elementary
geometry we see that if $(l,l',l'')$ are
$(\alpha,\beta)$-transverse, then there exists $k$ such that $l
\in L_k$ and $l',l'' \in L'_k$.  Thus the number of
$(\alpha,\beta)$-transverse joints can be bounded by
$$ \sum_{k=1}^K \# \{ p: p \in l, l', l'' \hbox{ for some } (\alpha,\beta)-\hbox{transverse } l \in L_k, l',l'' \in L'_k \}.$$
Next, observe from elementary geometry that if $l, l', l''$ are
$(\alpha,\beta)$-transverse in $L_k \cup L'_k$, then after
applying a dilation by $1/\beta$ in the $\omega_k$ direction, the
resulting lines become $c$-transverse for some $c \sim 1$ (here we
are using the hypothesis that $\alpha \sim 1$).  Applying
\eqref{n-bound} we conclude that \begin{eqnarray*} \begin{aligned}
\# \{ p: p \in l, l', l'' \hbox{ for some }&
(\alpha,\beta)-\hbox{transverse } l \in L_k, l',l'' \in L'_k \}\\&
= O_\eps( n^\eps ) (\# L_k)^{1/2} \# L'_k \end{aligned}
\end{eqnarray*} so it suffices to establish the bound
$$ \sum_{k=1}^K (\# L_k)^{1/2} \# L'_k \leq C n^{3/2} \beta^{-1/2}.$$
Now observe from transversality of $S$ and $S'$ that if $l$ has
direction in $S$ and $l'$ has direction in $S'$ then there are
at most $O(1)$ values of $k$ for which $l \in L_k$ and $l' \in
L_{k}'$.  This leads to the bound
$$ \sum_{k=1}^K \# L_k \# L'_k \leq C n^2.$$
On the other hand, observe that every line $l'$ belongs to at most
$O(1/\beta)$ families $L_{k}'$.  This leads to the bound
$$ \sum_{k=1}^K \# L'_k \leq C n \beta^{-1}.$$
The claim now follows from the Cauchy-Schwarz inequality.

Finally, we handle the case when $\alpha$ is very small, using the
bilinear rescaling argument from \cite{TVV}. Let $\tilde
\omega_1,\ldots,\tilde \omega_{\tilde K}$ be a maximal
$\alpha$-separated set of directions of the sphere, and for each
$\tilde \omega_k$ let $\tilde L_k$ be all the lines in $L$ which
make an angle of $O(\alpha)$ with $\tilde\omega_k$.  Observe that if
$(l,l',l'')$ are $(\alpha,\beta)$-transverse, then there exists
$k$ such that all of $l,l',l''$ lie in $\tilde L_k$. Thus we can
bound the total number of $(\alpha,\beta)$-joints in this case by
$$ \sum_{k=1}^{\tilde K} \# \{ p: p \in l, l', l'' \hbox{ for some } (\alpha,\beta)-\hbox{transverse } l, l',l'' \in \tilde L_k \}.$$
Next, observe from elementary geometry that if $l,l',l'' \in
\tilde L_k$ are $(\alpha,\beta)$-transverse, then if we dilate $l,
l', l''$ in the directions orthogonal to $\tilde \omega_k$ by
$1/\alpha$, then the resulting triple of lines becomes $(1,
\beta/\alpha)$-transverse. Since we have already established the
desired bound in the $\alpha \sim 1$ case, we conclude that
\begin{eqnarray*} \begin{aligned}
\# \{ p: p \in l, l', l'' \hbox{ for some }&
(\alpha,\beta)-\hbox{transverse } l, l',l'' \in \tilde L_k \}\\&
\leq O_\eps(n^\eps) (\# \tilde L_k)^{3/2+\eps}
(\beta/\alpha)^{-1/2-\eps}\end{aligned}\end{eqnarray*} and so it
will suffice to show that
$$ \sum_{k=1}^{\tilde K} (\# \tilde L_k)^{3/2+\eps} \leq C n^{3/2+\eps}.$$
Using the crude bound $\# \tilde L_k \leq n$, it suffices to show
that
$$ \sum_{k=1}^{\tilde K} \# \tilde L_k \leq C n.$$
But it is clear that each line $l \in L$ can belong to at most
$O(1)$ families $\tilde L_k$, and the claim follows.
\end{proof}

\begin{remark}
If one had the endpoint $q=d/(d-1)$ in Conjecture \ref{MLKC} then
one could remove the epsilon losses from the $n$ exponent, and
possibly also from the $\theta$ exponent as well.  The
deterioration of the bound as $\theta \to 0$ is closely related to
the reason that the multilinear Kakeya estimate is currently
unable to imply any corresponding linear Kakeya estimate.  Thus a
removal of this $\theta$-dependence in the joints estimate may
lead to a new \emph{linear} Kakeya estimate.
\end{remark}

\begin{remark} One can also phrase the joints problem for other families of
curves than lines, in the spirit of Section \ref{variable-sec}.  If one could
remove the loss of $\delta^{-\eps}$ in Theorem \ref{kakdiff}, one could obtain a result
similar to Theorem \ref{joints-thm} in this setting, but as Theorem \ref{kakdiff} stands
one would only obtain a rather unaesthetic result in which certain ``entropy numbers'' of
the joints are controlled.  We omit the details.
\end{remark}

\section{Appendix: A polynomial extrapolation lemma}

The main aim of this paper (the contents of Sections 3 and 4) is
to obtain monotonicity \textit{formulae} for spatial $L^{p}$-norms
of certain multilinear expressions. As we have seen, this can be
done quite explicitly for integer values of the exponent $p$, and
in such a way that the identities obtained \textit{make sense} at
least for non-integer $p$. The pay-off of having proved such
precise \textit{identities} for $p\in\N$ is that we may use a
density argument (e.g. using the
Weierstrass approximation theorem) to deduce that they must also
hold for $p\not\in\N$. This is very much analogous to the
classical result that a compactly supported probability
distribution is determined uniquely by its moments.

\begin{remark}
As we noted in Section 3, there is a satisfactory way of
avoiding this ``integer $p$ first" approach to the unperturbed
situation (Theorem \ref{lwmon}). This involves finding an
appropriate function of divergence form which differs from the
integrand in \eqref{train} by a manifestly non-negative quantity.  See \cite{BCCT}, \cite{CLL}.
In principle one could take a similar approach to Theorem
\ref{wowgen}, although as yet it seems quite unclear how to
directly exhibit an appropriate divergence term.
\end{remark}

\begin{lemma}\label{densitycor}
Suppose $f_{1},\hdots,f_{n}:\R^d\rightarrow \R$ are non-negative
bounded measurable functions for which the product $f_{1}\cdots
f_{n}$ is rapidly decreasing. Suppose that $G_{1},G_{2}:\R^{n}
\times\R^{d}\rightarrow\R$ are polynomial in their first variables
$p=(p_1,\hdots,p_n)$, with coefficients which are measurable and
of polynomial growth in their second. Then if the identity
$$
\int_{{\bf R}^{d}}G_{1}(p,x)f_{1}(x)^{p_{1}}\cdots f_n(x)^{p_n}dx=
\int_{{\bf R}^{d}}G_{2}(p,x)f_1(x)^{p_1}\hdots f_n(x)^{p_n}dx,
$$
holds for all $p\in \N^n$, then it holds for all $p\in
(0,\infty)^{n}$.
\end{lemma}


\begin{proof}
By linearity we may assume that $G_2=0$ and rename $G_1$ as $G$. Write
$$G(p,x) = \sum_{|\alpha| \leq N} p^{\alpha} w_{\alpha} (x)$$
where $w_{\alpha}$ is measurable and of polynomial growth.
 We may further
assume
that $\|f_j\|_{\infty} = 1$ for all $j$. Since the function
$$p \mapsto \int_{{\bf R}^{d}}G(p,x)f_{1}(x)^{p_{1}}\cdots f_n(x)^{p_n}dx$$
is an analytic function of each $p_j>0$, it suffices to prove the
result when $p_j \geq N$ for all $j$.
Let $\phi_p(t) := \phi_p (t_1, t_2, \dots , t_n) = t_1^{p_1} \dots t_n^{p_n}$
for $t \in [0,1]^n$. Note that $\phi_p (t) = 0$ if any $t_j = 0.$
Observe that we may write
$$G(p,x)f^p(x) = \sum_{|\alpha| \leq N} \tilde{w}_\alpha (x) \phi_p ^{(\alpha)}
(f(x)) $$
where $\tilde{w}_\alpha (x) := r_{\alpha} (x)
\prod_{j: \alpha_j \neq 0} f_j (x)$, and
$\phi_p^{(\alpha)}$ denotes differentiation of order $\alpha$ and where
$r_\alpha$ is of polynomial growth. In particular, if each $\alpha_j \geq 1$,
$\tilde{w}_\alpha \in L^1.$
By hypothesis, if $\phi$ is any polynomial which vanishes on the coordinate
axes,
\begin{equation}\label{last}
\int_{{\bf R}^d} \sum_{|\alpha| \leq N} \tilde{w}_\alpha (x)
\phi^{(\alpha)}(f(x)) dx = 0.
\end{equation}
We wish to show that the same continues to hold for $\phi$ replaced by
$\phi_p$ when each $p_j \geq N$. For such $\phi_p$ (which belong to the class
$ {\mathcal{C}} := \{\psi \in C^{N} ([0,1]^n)  : \psi(t)=0 \;
\mbox{whenever some} \; t_j =0\})$
we can approximate it to within any given $\eps$ by a polynomial $\phi$ of
class $\mathcal{C}$ in the norm $\| \psi \|_{*} :=
\max\{\| \psi^{(\alpha)}\|_{\infty}:\alpha_j \geq 1 \;
\forall j\mbox{ and }|\alpha|\leq N\}$.
So
$$ \int_{{\bf R}^d} \sum_{|\alpha| \leq N} \tilde{w}_\alpha (x)
\phi_p ^{(\alpha)}(f(x)) dx =
\int_{{\bf R}^d} \sum_{|\alpha| \leq N} \tilde{w}_\alpha (x)
[\phi_p ^{(\alpha)} - \phi^{(\alpha)}](f(x)) dx \; +0.$$
When $\alpha$ is such that each $\alpha_j \geq 1$ we can dominate
its contribution to the right hand side by
$\int_{{\bf R}^d} | \tilde{w}_{\alpha} (x)| dx \; \| \phi_p - \phi \|_{*}$
which is as small as we like since $\tilde{w}_{\alpha}$ is in $L^1$.

When some of the $\alpha_j$ are zero, say $\alpha_1, \dots , \alpha_k = 0$ and
$\alpha_{k+1}, \dots ,\alpha_n \neq 0$, we set $\tilde{\alpha} = (1,1, \dots,
1, \alpha_{k+1}, \dots, \alpha_n)$ and write
$$\phi_p ^{(\alpha)}(t) - \phi^{(\alpha)}(t) =
\int_0^{t_1} \dots \int_0^{t_k} [\phi_p ^{(\tilde{\alpha})} -
\phi^{(\tilde{\alpha})}](s_1, \dots s_k , t_{k+1}, \dots , t_n)
ds_1 \dots ds_k $$

Thus $[\phi_p ^{\alpha} - \phi ^{\alpha}](f(x)) \leq f_1(x) \dots f_k (x)
\|\phi_p -\phi\|_{*}$ and so for these $\alpha$
\begin{eqnarray*}
\begin{aligned}
\int_{{\bf R}^d}
|\tilde{w}_{\alpha} (x) | |[\phi_p ^{\alpha} - &\phi ^{\alpha}]
(f(x))| dx\\
& \leq \int_{{\bf R}^d} f_1(x) \dots f_k (x) f_{k+1}(x) \dots f_n(x)
|r_\alpha (x)| dx
\| \phi_p - \phi \|_{*}
\end{aligned}
\end{eqnarray*}
which is likewise as small as we like since $r_{\alpha}$ is of polynomial
growth and $f_1 \dots f_n$ is rapidly decreasing.

Thus formula \eqref{last} continues to hold for $\phi_p$ and we are finished.
\end{proof}

\end{document}